\documentclass{AIMS}
\usepackage{amssymb,amsmath}
\usepackage{float,bm}
  \usepackage{paralist}
  \usepackage{graphics} 
  \usepackage{epsfig} 
\usepackage{graphicx}  \usepackage{epstopdf}
 \usepackage[colorlinks=true]{hyperref}
\hypersetup{urlcolor=blue, citecolor=red}
\usepackage{mathrsfs}
\usepackage{tikz}
\usepackage{pgfplots}
\usetikzlibrary{positioning,shapes.geometric,arrows.meta,calc,fit}
\pgfplotsset{compat=1.18}

\definecolor{light_gray}{gray}{0.75}
\definecolor{dark_green}{rgb}{0.0, 0.6, 0.0}
\definecolor{royal_blue}{rgb}{0.0, 0.22, 0.66}
\definecolor{gold}{rgb}{0.8, 0.63, 0.21}
\definecolor{crimson}{rgb}{0.79, 0.0, 0.09}
\definecolor{amethyst}{rgb}{0.6, 0.4, 0.8}
\definecolor{copper}{rgb}{0.72, 0.45, 0.2}
\definecolor{brown}{rgb}{0.59, 0.29, 0.0}
\definecolor{forestgreen}{rgb}{0.133, 0.543, 0.113}
\definecolor{azure}{rgb}{0.0, 0.5, 1.0}
\definecolor{salmon}{rgb}{1.0, 0.55, 0.41}
\definecolor{sea_green}{rgb}{0.02, 0.52, 0.51}
\definecolor{teal_custom}{rgb}{0.0, 0.5, 0.5}

\newcommand{\CHANGEF}[1]{{ #1}}

\newcommand{\bb}{\boldsymbol b}
\newcommand{\bn}{\boldsymbol n}

\newcommand{\GD}{\Gamma^D}
\newcommand{\GN}{\Gamma^N}

\newcommand\calG{\mathscr{G}}
\newcommand\calT{\mathscr{T}}

\newcommand{\Pe}{\mathrm{Pe}}
\newcommand{\deltaSTD}{\delta_K^{\mathrm{std}}}
\newcommand{\deltaBND}{\delta_K^{\mathrm{bnd}}}
\newcommand{\deltaOPT}{\delta_K^{\mathrm{opt}}}
\newcommand{\deltaNN}{\delta_K^{\mathrm{NN}}}
\newcommand{\uSUPG}{u_h^\mathrm{SUPG}}
\newcommand{\uTAB}{u_h^\mathrm{Tabata}}
\newcommand{\uDIFF}{u^{\mathrm{diff}}_K}
\newcommand{\xiINT}{\xi^{\mathrm{int}}}
\newcommand{\xiBND}{\xi^{\mathrm{bnd}}}
\newcommand{\xiOPT}{\xi^{\mathrm{opt}}}
\newcommand{\Lmain}{{\mathrm L}_{\mathrm{main}}}
\newcommand{\Laux}{{\mathrm L}_{\mathrm{aux}}}
\newcommand{\xiMain}{\xi_{\mathrm{main}}}
\newcommand{\xiAux}{\xi_{\mathrm{aux}}}

\newtheorem{example}{Example}

\def\currentvolume{X}
 \def\currentissue{X}
  \def\currentyear{200X}
   \def\currentmonth{XX}
    \def\ppages{X--XX}
     \def\DOI{10.3934/xx.xx.xx.xx}

\title[Machine Learning for SUPG Parameters on Unstructured Meshes] 
{Application of Machine Learning to the Computation of Parameters in the 
SUPG Method for Unstructured and Anisotropic Meshes}

\author[P. Knobloch and M. Prakash]{}


\thanks{$^*$ Corresponding author: Petr Knobloch}

\begin{document}

\maketitle

\centerline{\scshape Petr Knobloch$^*$ and Manoj Prakash}
\medskip
{\footnotesize
\centerline{Department of Numerical Mathematics, Faculty of Mathematics and 
Physics}
\centerline{Charles University}
\centerline{Sokolovsk\'{a} 83, 186 75 Prague 8, Czech Republic}
\centerline{e-mails: knobloch@karlin.mff.cuni.cz, manoj@karlin.mff.cuni.cz}
}

\bigskip

 \centerline{(Communicated by the associate editor name)}

\begin{abstract}
Solution of convection-dominated problems using stabilized methods often
requires to specify stabilization parameters whose optimal choice is not known
but which considerably influence the quality of the approximate solution. In a
recent work of the authors, an approach for computing these parameters by
machine learning was proposed and applied to the streamline
upwind/Petrov--Galerkin (SUPG) method for convection--diffusion equations. In
the present paper, this approach is studied numerically for unstructured meshes
and extended to meshes obtained by anisotropic mesh adaptation.
\end{abstract}

\section{Introduction}
\label{sec:intro}

Many physical processes are influenced by convective and diffusive effects, and
the convective transport often dominates the diffusion by several orders of
magnitude. This leads to difficulties when such a physical process is simulated
numerically so that special discretization techniques have to be applied. In the
context of the finite element method, stabilization terms involving user-chosen
parameters are often introduced, see, e.g., \cite{RST08}. Since it is difficult
to find optimal values of these parameters, an approach based on machine
learning was recently proposed in \cite{KP26} and applied to the SUPG method,
leading to first promising numerical results. In this paper, we concentrate on
investigating the capabilities of this technique for unstructured meshes and
extend it to meshes obtained by anisotropic mesh adaptation.

As a model situation, we consider the convection--diffusion problem
\begin{equation}\label{steady-strong}
   -\varepsilon\,\Delta u + \bb\cdot\nabla u = f \quad\textrm{in}\;\Omega\,,
   \qquad
   u = u_b\quad\textrm{on}\;\GD\,,\qquad
   \varepsilon\,\frac{\partial u}{\partial\bn} = g\quad\textrm{on}\;\GN\,, 
\end{equation}
where $\Omega\subset{\mathbb R}^2$ is a bounded domain with a polygonal
Lipschitz--continuous boundary $\partial\Omega$ which is decomposed into a
Dirichlet part $\Gamma^D$ and a Neumann part $\Gamma^N$. We assume that
$\Gamma^D$ has a positive one-dimensional measure. The symbol $\bn$ denotes the 
outward-pointing unit normal vector to $\partial\Omega$. The data in
\eqref{steady-strong} have the following physical meanings: $\varepsilon>0$ is
the diffusivity which is assumed to be constant, $\bb\in
W^{1,\infty}(\Omega)^2$ is the convection field which is assumed to be
divergence-free, $f\in L^2(\Omega)$ describes sources and sinks of the quantity
$u$, and the functions $u_b\in H^{1/2}(\Gamma^D)$ and $g\in H^{-1/2}(\Gamma^N)$
specify the Dirichlet and Neumann boundary conditions, respectively. Moreover,
we assume that the inflow boundary is a part of the Dirichlet boundary, i.e.,
$\{\mathbf x\in\partial\Omega\,;\,\,(\bb\cdot\bn)(\mathbf x)<0\}\subset
\overline{\Gamma^D}$. Under these assumptions, the problem
\eqref{steady-strong} has a unique weak solution $u\in H^1(\Omega)$, see, e.g.,
\cite{GT01}.

In typical applications, convection dominates diffusion, which causes that the
solution of \eqref{steady-strong} usually contains so-called layers which are
narrow regions where the solution $u$ changes abruptly. As already mentioned, 
it is then difficult to solve \eqref{steady-strong} numerically since standard
methods often lead to spurious oscillations and unphysical values in the
approximate solution. In this paper, we concentrate on the finite element
method where many approaches for the numerical solution of
\eqref{steady-strong} were developed by modifying the standard Galerkin
discretization, see, e.g., the monographs \cite{RST08,BJK25}. Often, these
approaches are referred to as stabilized methods. One of the first stabilized
methods was the streamline upwind/Petrov--Galerkin (SUPG) method proposed in
\cite{HB79,BH82} which is still very popular due to its stability properties
and high accuracy away from layers. Nevertheless, like for many other linear
stabilized methods, the approximate solutions computed using the SUPG method
still contain unphysical values. However, it is important that, in contrast to
the Galerkin method, these unphysical values are restricted to layer regions
and bounded independently of $\varepsilon$. The SUPG method contains
stabilization parameters which determine the amount of the artificial diffusion
added by the stabilization term and hence influence the quality of the
approximate solution. Since an optimal choice of the stabilization parameters
is not known, an a posteriori parameter optimization framework was developed in
\cite{JKS11}. For model problems, this approach allows to obtain approximate
solutions of remarkable quality \cite{JK13,KLS19} but the method is
computationally quite expensive and its success depends on the difficult task
of choosing an appropriate target functional, see the studies in \cite{JKW23}.
Therefore, it was proposed in \cite{KP26} to compute the stabilization
parameters using machine learning, i.e., to approximate the functional
dependence between local features and the desired stabilization parameter using
a neural network. The costly parameter optimization approach is then used in
the training phase only. In addition, it is then much easier to define an
appropriate target functional since a reference solution to be approximated
using the SUPG method with optimized parameters is known. The local features
are determined by the data of \eqref{steady-strong}, some geometric information
and two cheaply computable finite element approximations of $u$ obtained by the
Tabata method \cite{Tabata77} and the SUPG method with a standard stabilization
parameter. Consequently, using a neural network allows to predict appropriate
stabilization parameters with reduced computational cost in comparison to
existing methods and with undiminished accuracy of the solution.

The use of machine learning to compute SUPG stabilization parameters has been
explored from several angles in recent years.
\CHANGEF{An early contribution is \emph{SPDE-Net} \cite{yadav2021spdenet},
where the stabilization parameter is predicted by an artificial neural network
trained either in a supervised manner or by minimizing a global or local $L^2$
error measure. This line of work was subsequently developed towards fully
residual-based training in \cite{MR4678304}, where the parameter is}
the output of a neural network whose loss is built from the residual of the
partial differential equation (PDE) itself. It is a physics-informed neural
network with no pre-computed optimized parameters and the network learns the
stabilization parameter by being penalized for producing SUPG solutions that
violate the strong form of the equation. This idea was extended in
\cite{yadav2024convstabnet}, where a convolutional architecture was used to
incorporate spatial information from neighbouring cells. In
\cite{yadav2023qpde}, the network was replaced by a quantum neural network. The
choice of the architecture and training hyperparameters within this paradigm
was studied in \cite{joshi2021choice}.
\CHANGEF{A related supervised approach was proposed in \cite{tassi2023}, where
the training data are generated by solving optimization problems that minimize
the distance between the numerical and the exact solution for varying problem
data and finite element settings, and the trained network is then used online
to predict the stabilization parameter.}

\CHANGEF{The approach followed in the present paper differs from these works in
several respects. First, the parameter is predicted element-wise, so that the
amount of artificial diffusion adapts to the local behaviour of the solution
rather than to the global setting of the problem. Second, the input features
include information extracted from the Tabata solution \cite{Tabata77}, which
satisfies the discrete maximum principle and therefore provides a cheap
indicator of both layer regions and spurious oscillations; the same solution is
used as a reference in a post-processing step. Third, mesh cells intersecting
the outflow Dirichlet boundary are treated by a separate model, since the
optimal parameter behaves differently there. Finally, and in contrast to all
the works mentioned above, we consider unstructured and anisotropically adapted
meshes, in which the elements are strongly stretched.}

A separate line of research applies neural networks to convection-diffusion
problems by replacing the discretization itself rather than augmenting it, for
instance through  physics-informed neural networks
\cite{anandh2025improving,hasan2025mesh} or specialized deep-learning\linebreak
schemes \cite{beguinet2023deep}. These methods solve the PDE directly with the
neural network and are therefore complementary to our approach. We refer to
\cite{KP26} for a review of further publications related to the present paper.

The plan of the paper is as follows. In the next section, we formulate the SUPG
method, mention some of its properties and present formulas for the
stabilization parameter. Section~\ref{sec:ml_iso} reviews the main ideas of our
approach of computing the stabilization parameters using machine learning. This
approach is then studied numerically in Section~\ref{sec:hemker} for the Hemker
problem and a sequence of unstructured meshes. Section~\ref{sec:ml_aniso}
extends the ideas presented in Section~\ref{sec:ml_iso} to meshes obtained by
anisotropic mesh adaptation. The properties of this new approach are
illustrated by numerical results in 
Section~\ref{sec:numerics_aniso}
\CHANGEF{and the influence of some components of the method is
examined in Section~\ref{sec:ablation}. Finally, Section~\ref{sec:conclusions}
summarizes the results and outlines directions for future work.}

\section{SUPG method}
\label{sec:supg}

To define the SUPG method, we introduce a triangulation $\calT_h$ decomposing
$\overline\Omega$ into a finite number of triangles possessing the usual
compatibility properties. The discretization parameter $h$ represents the
maximum of the diameters of the triangles in $\calT_h$. We define the finite
element space
\begin{equation*}
   W_h=\{v_h\in C(\overline\Omega)\,;\,\,
   v_h|_K^{}\in P_1(K)\,\,\,\forall\,\,K\in\calT_h\}\,,
\end{equation*}
where $P_1(K)$ is the space of linear functions on $K$, and denote
$$
   V_h=\{v_h\in W_h\,;\,\,v_h=0\,\,\,\,\mbox{on}\,\,\,\Gamma^D\}\,.
$$
Moreover, we introduce a function $u_{bh}\in W_h$ such that its trace
on $\Gamma^D$ approximates the boundary condition $u_b$. Then the SUPG
discretization of \eqref{steady-strong} reads: Find $u_h\in W_h$ such that
$u_h=u_{bh}$ on $\Gamma^D$ and
\begin{align}
   &\varepsilon\,(\nabla u_h,\nabla v_h)+(\bb \cdot\nabla u_h,v_h)
   +\sum_{K\in\calT_h}\delta_K\,
   (-\varepsilon\,\Delta u_h + \bb\cdot\nabla u_h - f,\bb\cdot\nabla v_h)_K
   \nonumber\\
   &\hspace*{55mm}=(f,v_h)+\langle g,v_h\rangle_{\Gamma^N}\quad
   \forall\,\,v_h\in V_h\,,\label{eq:supg}
\end{align}
where $\delta_K\ge0$ is a stabilization parameter, constant in each mesh cell
$K\in\calT_h$, $(\cdot,\cdot)$ denotes the inner product in $L^2(\Omega)$ or
$L^2(\Omega)^2$, $(\cdot,\cdot)_K$ denotes the inner product in $L^2(K)$, and
$\langle\cdot,\cdot\rangle_{\Gamma^N}$ is the duality pairing between
$H^{-1/2}(\Gamma^N)$ and $H^{1/2}(\Gamma^N)$. We consider the term $\Delta u_h$
in \eqref{eq:supg} although it vanishes on any $K\in\calT_h$ to stress the
residual-based character of the stabilization which assures the consistency of
the method. We observe that the SUPG stabilization term introduces artificial
diffusion in the streamline direction, which enhances the stability of the
method. More precisely, the SUPG method is stable with respect to the norm
\begin{displaymath}
   |||v|||=\left(\varepsilon\,|v|_{1,\Omega}^2 
     +\|\delta_h^{1/2}\,\bb\cdot\nabla v\|_{0,\Omega}^2
     +\frac12\,\|(\bb\cdot\bn)^{1/2}\,v\|_{0,\Gamma^N}^2\right)^{1/2}.
\end{displaymath}

To guarantee a convergence of the SUPG method, the stabilization parameters
$\delta_K$ have to satisfy some bounds dependent on the mesh and the
data of \eqref{steady-strong}.
\CHANGEF{As an example, we present the setting considered in \cite{KP26} where
it was required}
that, for any $K\in\calT_h$, the value $\delta_K$ satisfies
\begin{align}
   &\delta_K\le\alpha\,\frac{h_K^2}
    {\varepsilon+h_K\,\|\bb\|_{0,\infty,K}^{}}\,,\label{upper}\\[1mm]
   &\delta_K\ge\frac1\alpha\,\frac{h_K}{\|\bb\|_{0,\infty,K}^{}}\qquad
    \mbox{if}\quad\frac{h_K\,\|\bb\|_{0,\infty,K}^{}}\varepsilon>\alpha\,,
    \label{lower}
\end{align}
where $h_K=\mathrm{diam}(K)$ and $\alpha>4/3$ is a fixed constant independent
of $K$ and~$h$. Moreover,
\CHANGEF{it was assumed in \cite{KP26}}
that $\calT_h$ belongs to a family of shape-regular triangulations in the sense
that
\begin{equation}\label{eq:shape_reg}
    \frac{h_K}{\varrho_K}\le\sigma\qquad\forall\,\,K\in\calT_h\,,
\end{equation}
where $\varrho_K$ is the diameter of the largest circle inscribed in $K$ and 
$\sigma>1$ is a constant independent of $K$ and $h$. Finally, assume that
$u_b\in C(\overline{\Gamma^D})$ and that $u_{bh}$ coincides with $u_b$ at the
vertices of $\calT_h$ contained in $\overline{\Gamma^D}$. Under the above
assumptions, it was proved in \cite{KP26} that, if the weak solution of
\eqref{steady-strong} satisfies $u\in H^2(\Omega)$, then the solution of the
SUPG discretization \eqref{eq:supg} fulfills the error estimate
\begin{equation*}
   |||u-u_h|||\le\alpha^{1/2}\,C\,
   (\varepsilon+\|\bb\|_{0,\infty,\Omega}^{}\,h)^{1/2}\,h\,|u|_{2,\Omega}^{}\,,
\end{equation*}
where the constant $C$ depends only on $\sigma$ from \eqref{eq:shape_reg}. 
\CHANGEF{For anisotropic meshes (which do not satisfy the assumption
\eqref{eq:shape_reg}), the analysis as well as the resulting estimates are much
more complicated and only few results can be found in the literature, see, e.g.,
\cite{MPP03}.}

There are various strategies for choosing the stabilization parameters
$\delta_K$, see, e.g., the discussions in \cite{JK07} and \cite{KP26}. Often,
one uses
\begin{equation}\label{eq:delta_std}
   \deltaSTD=\frac{h_K^\parallel}{2\,\|{\bb}\|_{0,\infty,K}^{}}
    \left(\coth \Pe_K-\frac1{\Pe_K}\right)\quad\mbox{with}\quad 
    \Pe_K= \frac{h_K^\parallel\,\|{\bb}\|_{0,\infty,K}^{}}{2\,\varepsilon}\,,
\end{equation}
where 
\begin{equation}\label{eq:hKparallel}
   \CHANGEF{h_K^\parallel:=\max\{|{\mathbf x}-{\mathbf y}|\,;\,\,
   {\mathbf x},{\mathbf y}\in K,\,
   {\mathbf x}-{\mathbf y}=\alpha\,\bb({\mathbf x}_K),\,\alpha\in{\mathbb R}\}}
\end{equation}
is the diameter of $K$ in the direction of the convection $\bb$ evaluated at
the barycentre
\CHANGEF{${\mathbf x}_K$}
of $K$. If the convection vanishes at
\CHANGEF{${\mathbf x}_K$},
we set $h_K^\parallel:=\mathrm{diam}(K)$. The symbol $\Pe_K$
denotes the local P\'eclet number that determines whether the problem is
locally (i.e., within a particular mesh cell) convection-dominated or
diffusion-dominated. The formula \eqref{eq:delta_std} is a generalization of a
formula derived in the one--dimensional case \cite{CGMZ76} to obtain a nodally
exact solution for constant data and a uniform division of $\Omega$. We will
call $\delta_K$ defined by \eqref{eq:delta_std} standard stabilization
parameter. Note that $\delta_K$ defined by \eqref{eq:delta_std} satisfies the
upper bound \eqref{upper} with $\alpha=1$ and,
\CHANGEF{under the shape-regularity assumption \eqref{eq:shape_reg}, it also
satisfies}
the lower bound \eqref{lower} with $\alpha=5\,\sigma$.

For outflow boundary layers, an improved definition of $\delta_K$ was
introduced in \cite{Kno09b} and further explored in \cite{Kno08}. To formulate
this approach, let us introduce the outflow Dirichlet boundary
\begin{displaymath}
   \Gamma^{\mathrm{out}}=\overline{
   \{{\mathbf x}\in\Gamma^D\,;\,\,(\bb\cdot\bn)({\mathbf x})>0\}}\,,
\end{displaymath}
and set
\begin{displaymath}
   G_h=\bigcup_{K\in\calG_h}K\qquad\mbox{where}\qquad 
   {\calG}_h=\{K\in\calT_h\,;\,\,K\cap\Gamma^{\mathrm{out}}\neq\emptyset\}\,.
\end{displaymath}
It was shown in \cite{Kno09b} that it is possible to construct a piecewise constant function $\delta_0$ on $G_h$ satisfying
\begin{equation*}
   \int_{G_h}\,v_h+\delta_0\,\bb\cdot\nabla v_h\,{\rm d}{\mathbf x}=0\qquad
   \forall\,\,v_h\in V_h\,.
\end{equation*}
Then, by analogy to \eqref{eq:delta_std}, we define the parameter $\delta_K$,
on any element $K\in{\calG}_h$, by
\begin{equation}\label{eq:delta0_bnd}
   \deltaBND=\delta_0|_K^{}\left(\coth \Pe_K-\frac1{\Pe_K}\right).
\end{equation}
This choice suppresses the influence of the Dirichlet boundary condition at the
outflow boundary $\Gamma^{\mathrm{out}}$ on the values of $u_h$ at interior
nodes, which considerably reduces unphysical values of $u_h$ along
$\Gamma^{\mathrm{out}}$. On elements 
$K\in{{\calT}_h}\setminus{\calG}_h$, we set $\delta_K=\deltaSTD$.

\section{Computation of $\delta_K$ using machine learning}
\label{sec:ml_iso}

As already mentioned in the introduction, the accuracy of the SUPG method can be
significantly improved by computing the stabilization parameters $\delta_K$
using the adaptive a posteriori optimization approach proposed in \cite{JKS11}.
Let us explain the basic idea of this approach. We denote by $\delta_h$ a
piecewise constant function such that $\delta_h|_K^{}=\delta_K$ for any
$K\in\calT_h$. All permissible functions $\delta_h$ form a set $D_h\subset
L^\infty(\Omega)$. This set can be determined, e.g., by requiring the validity
of the bounds \eqref{upper} and \eqref{lower} with some appropriate
constant~$\alpha$. Then, for any $\delta_h\in D_h$, there is a uniquely
determined solution $u_h=u_h(\delta_h)$ of the SUPG discretization
\eqref{eq:supg}. Furthermore, we introduce a target functional $I_h$ which
measures the quality of the solution $u_h$. Assuming that a smaller value of
$I_h(u_h)$ means a better solution $u_h$, the aim is to find $\delta_h\in D_h$
such that $I_h(u_h(\delta_h))$ is small, i.e., to minimize
$\Phi(\delta_h):=I_h(u_h(\delta_h))$ over $D_h$. This is a constrained
nonlinear optimization problem.

The success of the described parameter optimization approach depends on the
used minimization algorithm and on the considered target functional $I_h$ which
has to correspond to the solution of \eqref{steady-strong}, see, e.g.,
\cite{JKS11,JK13,KLS19,JKW23}. Consequently, in general, it may be very
difficult to construct the target functional in an appropriate way. Therefore,
it was proposed in \cite{KP26} to use the parameter optimization algorithm for
learning a neural network designed to compute the stabilization parameters. In
this case, the training is based on a set of selected data defining the problem
\eqref{steady-strong}.
For each problem
\eqref{steady-strong} used in the training phase, an accurate approximation of
the solution $u$ can be obtained. Either it is directly prescribed together
with the data or it is computed by means of an accurate numerical method
satisfying the discrete maximum principle (DMP). Such methods are always
nonlinear, cf., e.g., \cite{BJK23,BJK25}, and hence the computation of an
accurate approximation of $u$ is costly. Therefore, it is performed in the
training phase only. Having an accurate approximation $\tilde u_h$ of~$u$,
which can be computed, e.g., by some of the algebraic flux correction schemes
from \cite[Chapter 10]{BJK25}, a target functional $I_h$ for computing the
optimized SUPG stabilization parameters $\deltaOPT$ is much easier to
construct. One can simply define it as a measure of the distance of the SUPG
solution to the accurate approximation $\tilde u_h$, see \cite{KP26}. In our
approach, the initial guess for $\deltaOPT$ is given by \eqref{eq:delta_std}
with the modification \eqref{eq:delta0_bnd} along outflow Dirichlet boundaries.

To design a neural network for computing the SUPG stabilization parameter
$\delta_K$ for a given $K\in\calT_h$, some cheaply computable local quantities
which are in correlation with the optimized  stabilization parameter
$\deltaOPT$ are needed. In \cite{KP26} we used local quantities obtained from
two cheaply computable finite element solutions. One of them was the solution
$\uSUPG$ of the SUPG method \eqref{eq:supg} with the standard stabilization
parameter~\eqref{eq:delta_std}. The other one was the solution $\uTAB$ of the
Tabata upwind scheme~\cite{Tabata77}. The Tabata method was the first upwind
finite element discretization satisfying the DMP, nevertheless, according to
the numerical studies in \cite{BJK25}, it is more accurate than other linear
upwind finite element methods satisfying the DMP developed later. We refer to
\cite[Section 8.4.1.1]{BJK25} for details on the definition of Tabata's method
and its properties.

Given $K\in\calT_h$, we compute the quantities
\begin{equation}\label{eq:u_quantities}
   \|\nabla\uSUPG\|_{L^\infty(K)}^{}\qquad\mbox{and}\qquad
   \|\uSUPG-\uTAB\|_{L^\infty(\Delta_K)}^{}\,,
\end{equation}
where
$$
   \Delta_K=\cup\{\tilde K\in\calT_h\,;\,\,\tilde K\cap K\neq\emptyset\}\,.
$$
The former quantity in \eqref{eq:u_quantities} helps to find layer regions
whereas the latter one serves for identifying spurious oscillations.
Furthermore, as observed in the preceding section, $\deltaOPT$ should depend on
$\varepsilon$ and $\|\bb\|_{L^\infty(K)}^{}$. Finally, $\deltaOPT$ should also
depend on characteristic dimensions of $K$, see, e.g., \cite{MPP03,CS07,DXY18}.
To this end, we consider the diameter $h_K^\parallel$ of~$K$ in the direction
of $\bb$ introduced
\CHANGEF{in \eqref{eq:hKparallel}}
and the diameter $h_K^\perp$ of~$K$ in the direction
\CHANGEF{$\bb^\perp$}
orthogonal to the convection $\bb$ evaluated at the barycentre
\CHANGEF{${\mathbf x}_K$} of~$K$.
\CHANGEF{The value $h_K^\perp$ is defined by \eqref{eq:hKparallel} with 
$\bb({\mathbf x}_K)$ replaced by $\bb^\perp({\mathbf x}_K)$.}
Again, $h_K^\perp:=\mathrm{diam}(K)$ if $\bb$ vanishes at
\CHANGEF{${\mathbf x}_K$}.
However, as it was shown in \cite{KP26}, scaling properties of the
stabilization parameter indicate that $\deltaOPT$ cannot depend on these six
quantities in an arbitrary way. A permissible type of dependence is
\begin{align}
   &\deltaOPT=\frac{(h_K^\parallel)^2}\varepsilon\,\xi\left(
   \frac{h_K^\parallel\,\|\bb\|_{L^\infty(K)}^{}}\varepsilon,
   \frac{h_K^\perp}{h_K^\parallel},
   \frac{h_K^\parallel\,\|\nabla\uSUPG\|_{L^\infty(K)}^{}}{\uDIFF},\right.
   \nonumber\\
   &\left.\hspace*{60mm}
   \frac{\|\uSUPG-\uTAB\|_{L^\infty(\Delta_K)}^{}}{\uDIFF}\right),
   \label{eq:delta_general}
\end{align}
where $\uDIFF$ is the difference between the maximum and minimum of the Tabata
solution $\uTAB$ in a neighbourhood of $K$. Here, for simplicity, we define
this scaling factor globally, i.e.,
$$
   \uDIFF=\max_{\overline\Omega}\,\uTAB-\min_{\overline\Omega}\,\uTAB\,.
$$

It turned out that the dependence \eqref{eq:delta_general} is appropriate for
mesh cells $K\in\calT_h\setminus\calG_h$, i.e., those ones which do not
intersect the outflow Dirichlet boundary. Thus, in particular,
\eqref{eq:delta_general} is used for interior mesh cells. In \cite{KP26}, we
considered a particular form of the dependence \eqref{eq:delta_general} given by
\begin{equation}\label{eq:delta_int}
   \deltaOPT=\deltaSTD\,\xiINT\left(
   \frac{h_K^\perp}{h_K^\parallel},
   \frac{h_K^\parallel\,\|\nabla\uSUPG\|_{L^\infty(K)}^{}}{\uDIFF},
   \frac{\|\uSUPG-\uTAB\|_{L^\infty(\Delta_K)}^{}}{\uDIFF}\right).
\end{equation}
For mesh cells $K\in\calG_h$ intersecting the outflow Dirichlet boundary, the
following modification of \eqref{eq:delta_int} was considered in \cite{KP26}:
\begin{equation}\label{eq:delta_bnd}
   \deltaOPT=\deltaBND\,\xiBND\left(
   \frac{h_K^\perp}{h_K^\parallel},
   \frac{h_K^\parallel\,\|\nabla\uSUPG\|_{L^\infty(K)}^{}}{\uDIFF},
   \frac{\|\uSUPG-\uTAB\|_{L^\infty(\widetilde\Delta_K)}^{}}{\uDIFF},
   \angle\bb_K\right),
\end{equation}
where 
$$
   \widetilde\Delta_K=\cup\{\tilde K\in\calT_h\,;\,\,
                      \tilde K\cap K\cap \Gamma^{\mathrm{out}}\neq\emptyset\}
$$
and
\CHANGEF{$\angle\bb_K:=\bb_K\cdot\bn_K$}
is a measure for the angle between the convection $\bb$ and
the normal vector to the boundary $\Gamma^{\mathrm{out}}$.
\CHANGEF{If $K\cap\Gamma^{\mathrm{out}}$ is an edge $E$ with the midpoint
${\mathbf x}_E$, then $\bb_K=\bb({\mathbf x}_E)/|\bb({\mathbf x}_E)|$ and
$\bn_K=\bn|_E^{}$. If $K\cap\Gamma^{\mathrm{out}}$ is a vertex ${\mathbf x}$,
then $\bb_K=\bb({\mathbf x})/|\bb({\mathbf x})|$ and $\bn_K$ is (an
approximation of) $\bn({\mathbf x})$.}

The idea of \cite{KP26} was to construct the nonlinear mappings $\xiINT$ and
$\xiBND$ in \eqref{eq:delta_int} and \eqref{eq:delta_bnd}, respectively, as
neural networks. Thus, the aim was to learn the neural networks the functional
relationships \eqref{eq:delta_int} and \eqref{eq:delta_bnd} between the
considered local features and the optimized stabilization parameter
$\deltaOPT$. Apart from the different number of input values, the architectures of the models and the training processes were
identical for both neural networks. Note that considering the functional
dependences \eqref{eq:delta_int} and \eqref{eq:delta_bnd} instead of general
dependences on the local quantities $\varepsilon$, $\|\bb\|_{L^\infty(K)}^{}$,
$h_K^\parallel$, $h_K^\perp$, $\|\nabla\uSUPG\|_{L^\infty(K)}^{}$,
$\|\uSUPG-\uTAB\|_{L^\infty(\Delta_K)}^{}$, and $\angle\bb_K$ significantly
simplifies the training and makes the learning process more stable and
efficient. By considering the relationships \eqref{eq:delta_int} and
\eqref{eq:delta_bnd}, the model can infer correct stabilization behavior even
in regimes not directly represented in the training data.

\begin{figure}[t!]
\centering
\scalebox{0.7}{%
\begin{tikzpicture}[
    font=\sffamily\small,
    >={Stealth[length=2mm]},
    neuron/.style={circle, draw, minimum size=5.5mm, inner sep=0pt, line width=0.4pt},
    inp/.style   ={neuron, fill=blue!15},
    hid/.style   ={neuron, fill=purple!20},
    gW/.style    ={neuron, fill=orange!20},   
    gN/.style    ={neuron, fill=orange!55},   
    mainN/.style ={neuron, fill=green!70},
    auxN/.style  ={neuron, fill=red!70},
    outM/.style  ={neuron, fill=green!75},
    outA/.style  ={neuron, fill=red!75},
    conn/.style  ={-, black!22, line width=0.15pt},
    gatebox/.style ={draw=orange!70, dashed, rounded corners=3pt, thick,
                     inner sep=5pt, line width=0.5pt},
    mainbox/.style ={draw=red!60,    dashed, rounded corners=3pt, thick,
                     inner sep=5pt, line width=0.5pt},
    auxbox/.style  ={draw=purple!60, dashed, rounded corners=3pt, thick,
                     inner sep=5pt, line width=0.5pt},
    ttl/.style     ={align=center, font=\small\bfseries},
    subttl/.style  ={align=center, font=\footnotesize\itshape, text=orange!60!black},
    cnt/.style     ={align=center, font=\footnotesize, text=black!70}
]

\def\xI{0}
\def\xH{2.0}
\def\xGa{3.8}  \def\xGb{5.4}  \def\xGc{7.0}    
\def\xGd{8.8}  \def\xGe{10.4} \def\xGf{12.0}   
\def\xBr{13.9}
\def\xO{16.0}

\foreach \i/\y in {1/1.00, 2/0.2} \node[inp] (I\i) at (\xI,\y) {};
\node at (\xI,-0.6) {$\vdots$};
\node[inp] (I3) at (\xI,-1.4) {};

\foreach \i/\y in {1/1.8, 2/1.0, 3/0.2, 4/-0.6} \node[hid] (H\i) at (\xH,\y) {};
\node at (\xH,-1.4) {$\vdots$};
\node[hid] (H5) at (\xH,-2.2) {};

\foreach \i/\y in {1/1.8, 2/1.0, 3/0.2, 4/-0.6} \node[gW] (Ga\i) at (\xGa,\y) {};
\node at (\xGa,-1.4) {$\vdots$};
\node[gW] (Ga5) at (\xGa,-2.2) {};

\foreach \i/\y in {1/1.0, 2/0.2, 3/-0.6} \node[gN] (Gb\i) at (\xGb,\y) {};
\node at (\xGb,-1.4) {$\vdots$};
\node[gN] (Gb4) at (\xGb,-2.2) {};

\foreach \i/\y in {1/1.8, 2/1.0, 3/0.2, 4/-0.6} \node[gW] (Gc\i) at (\xGc,\y) {};
\node at (\xGc,-1.4) {$\vdots$};
\node[gW] (Gc5) at (\xGc,-2.2) {};

\foreach \i/\y in {1/1.8, 2/1.0, 3/0.2, 4/-0.6} \node[gW] (Gd\i) at (\xGd,\y) {};
\node at (\xGd,-1.4) {$\vdots$};
\node[gW] (Gd5) at (\xGd,-2.2) {};

\foreach \i/\y in {1/1.0, 2/0.2, 3/-0.6} \node[gN] (Ge\i) at (\xGe,\y) {};
\node at (\xGe,-1.4) {$\vdots$};
\node[gN] (Ge4) at (\xGe,-2.2) {};

\foreach \i/\y in {1/1.8, 2/1.0, 3/0.2, 4/-0.6} \node[gW] (Gf\i) at (\xGf,\y) {};
\node at (\xGf,-1.4) {$\vdots$};
\node[gW] (Gf5) at (\xGf,-2.2) {};

\foreach \i/\y in {1/2.6, 2/1.8, 3/1.0} \node[mainN] (M\i) at (\xBr,\y) {};
\node at (\xBr,0.2) {$\vdots$};
\node[mainN] (M4) at (\xBr,-0.6) {};

\foreach \i/\y in {1/-2.0, 2/-2.8, 3/-3.6} \node[auxN] (X\i) at (\xBr,\y) {};
\node at (\xBr,-4.4) {$\vdots$};
\node[auxN] (X4) at (\xBr,-5.2) {};

\node[outM] (OM) at (\xO, 1.0)  {};
\node[outA] (OA) at (\xO,-3.6) {};

\foreach \a in {1,...,3} \foreach \b in {1,...,5} \draw[conn] (I\a) -- (H\b);
\foreach \a in {1,...,5} \foreach \b in {1,...,5} \draw[conn] (H\a) -- (Ga\b);
\foreach \a in {1,...,5} \foreach \b in {1,...,4} \draw[conn] (Ga\a) -- (Gb\b);
\foreach \a in {1,...,4} \foreach \b in {1,...,5} \draw[conn] (Gb\a) -- (Gc\b);
\foreach \a in {1,...,5} \foreach \b in {1,...,5} \draw[conn] (Gc\a) -- (Gd\b);
\foreach \a in {1,...,5} \foreach \b in {1,...,4} \draw[conn] (Gd\a) -- (Ge\b);
\foreach \a in {1,...,4} \foreach \b in {1,...,5} \draw[conn] (Ge\a) -- (Gf\b);
\foreach \a in {1,...,5} \foreach \b in {1,...,4} {
    \draw[conn] (Gf\a) -- (M\b);
    \draw[conn] (Gf\a) -- (X\b);
}
\foreach \b in {1,...,4} \draw[conn] (M\b) -- (OM);
\foreach \b in {1,...,4} \draw[conn] (X\b) -- (OA);

\node[gatebox, fit=(Ga1)(Gc5), label={[subttl]above:Attention Gate 1}] (G1box) {};
\node[gatebox, fit=(Gd1)(Gf5), label={[subttl]above:Attention Gate 2}] (G2box) {};
\node[mainbox, fit=(M1)(M4)(OM)] (Mbox) {};
\node[auxbox,  fit=(X1)(X4)(OA)] (Abox) {};

\node[ttl] at (\xI, 2.1) {Input};
\node[ttl] at (\xH, 2.9) {Hidden};
\node[ttl] at (\xBr, 3.4) {Main};
\node[ttl] at (\xBr,-6.0) {Auxiliary};
\node[ttl] at (\xO, 1.8) {Output};
\node[ttl] at (\xO,-2.8) {Output};

\node[cnt] at (\xI, -2.2) {3 or 4};
\node[cnt] at (\xH, -3.0) {256};
\node[cnt] at (\xGa,-3.0) {256};
\node[cnt] at (\xGb,-3.0) {64};
\node[cnt] at (\xGc,-3.0) {256};
\node[cnt] at (\xGd,-3.0) {256};
\node[cnt] at (\xGe,-3.0) {64};
\node[cnt] at (\xGf,-3.0) {256};
\node[cnt] at (\xBr, 3.0) {128};
\node[cnt] at (\xBr,-5.6) {128};
\node[cnt] at (\xO, 0.2)  {1};
\node[cnt] at (\xO,-4.4)  {1};

\end{tikzpicture}
}
\caption{Architecture of the neural network (MLP).}\label{fig:architecture}
\end{figure}
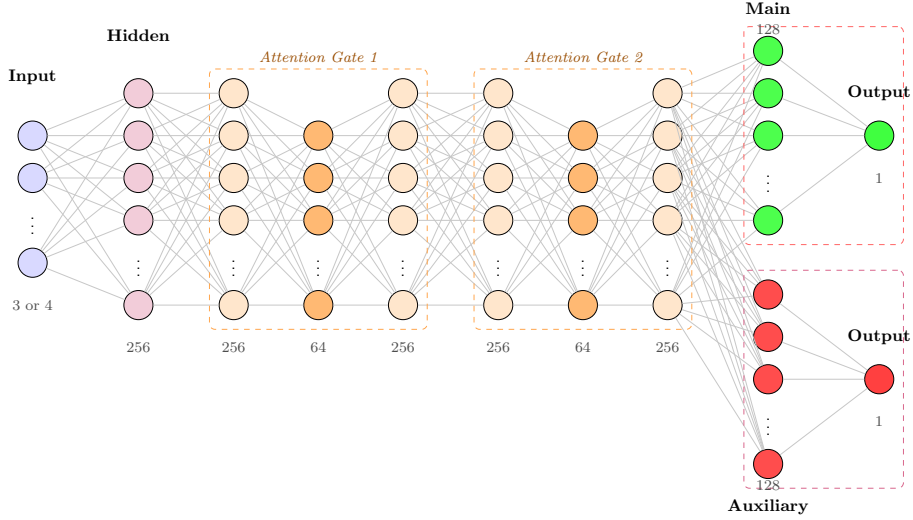

The neural network model considered in \cite{KP26} was developed using the
PyTorch library which provides tools that make it easier to build and train
deep learning models. The training data were generated from a diverse set of
examples that had either a boundary layer or an interior layer or both. The used
neural network architecture is depicted in Fig.~\ref{fig:architecture}.
The key component of the model is an encoder that transforms the input feature
vector of dimension 3 or 4 into a 256-dimensional vector. The encoder includes a
fully connected layer, followed by layer normalization to keep training stable,
and a ReLU activation to introduce nonlinearity. To make the model more robust
while avoiding too much regularization, it uses a custom AdaptiveDropout layer.
This layer changes the dropout rate based on the standard deviation of each
feature. The 256-dimensional vector from the encoder is passed through two
sequential AttentionGate modules. Each module includes a bottleneck network
($256\to64\to256$) that uses SiLU and Sigmoid activations. These gates help the
network focus more on important features and ignore less important ones at
different layers. All the weights are initialized using Kaiming normal
initialization (with a negative slope of 0.1 to maintain compatibility with
potential LeakyReLU variants), and biases are set to zero.

To improve information flow during training, the network uses deep supervision
by adding two output branches. The primary output head is a two-layer
multilayer perceptron (MLP) with dimensions $256\to128\to1$ with ReLU
activation. In parallel, an auxiliary head with an identical structure provides
an additional supervisory signal during training, which helps the model to train
more smoothly and keeps the internal feature representations more stable and
useful.
\CHANGEF{For each mesh cell $K$, the network does not predict the stabilization
parameter $\delta_K^{\mathrm{NN}}$ itself but the dimensionless multiplier of
\eqref{eq:delta_int} and \eqref{eq:delta_bnd}. For cells 
$K\in\calT_h\setminus\calG_h$ the raw output of the MLP is denoted by
$\xi_K^{\mathrm{int}}$, the training target is
\[
   \xi_K^{\mathrm{opt},\mathrm{int}}
   =\frac{\delta_K^{\mathrm{opt}}}{\delta_K^{\mathrm{std}}}\,,
   \qquad\text{and}\qquad
   \delta_K^{\mathrm{NN}}=\delta_K^{\mathrm{std}}\,\xi_K^{\mathrm{int}}\,.
\]
For cells $K\in\calG_h$ intersecting the outflow Dirichlet boundary, the raw
output of the network is denoted by $\xi_K^{\mathrm{bnd}}$, the training target is
\[
   \xi_K^{\mathrm{opt},\mathrm{bnd}}
   =\frac{\delta_K^{\mathrm{opt}}}{\delta_K^{\mathrm{bnd}}}\,,
   \qquad\text{and}\qquad
   \delta_K^{\mathrm{NN}}=\delta_K^{\mathrm{bnd}}\,\xi_K^{\mathrm{bnd}}\,.
\]
No further normalization or inverse transformation is applied, so that the
multiplier predicted by the network is the only learned quantity in both
cases.}

The loss functional for training is defined as a linear combination of two
losses, one from the main output and another one from the auxiliary output.
\CHANGEF{Denoting by $\xiMain$ and $\xiAux$ the multipliers predicted by
the main and the auxiliary head, respectively, and by $\xiOPT$ the target
multiplier, it is given by}
\begin{equation*}
   \CHANGEF{\mathrm{Loss}(\xiMain,\xiAux)
   =0.7\,\Lmain(\xiMain-\xiOPT)+0.3\,\Laux(\xiAux-\xiOPT)\,,}
\end{equation*}
where
\begin{equation*}
   \Lmain(a) = \begin{cases}
               \frac{1}{2}\,a^2 & \text{if}\,\,\, |a| \le 1\,, \\
               |a| - \frac{1}{2} & \text{otherwise}\,,
               \end{cases}\qquad
   \Laux(a) = a^2\,.
\end{equation*}
\CHANGEF{The auxiliary head is used during training only; at inference the
prediction is read from the main head, i.e., $\xi_K^{\mathrm{int}}$ and
$\xi_K^{\mathrm{bnd}}$ represent $\xiMain$.}

\CHANGEF{Since $\delta_K^{\mathrm{std}}\ge 0$ and $\delta_K^{\mathrm{bnd}}\ge
0$, the non-negativity $\delta_K^{\mathrm{NN}}\ge 0$ is inherited from the
non-negativity of the predicted multiplier and no clipping is applied; in
particular, the bounds \eqref{upper}, \eqref{lower} are not enforced
explicitly. If $\delta_K^{\mathrm{std}}=0$ or $\delta_K^{\mathrm{bnd}}=0$, the
target multiplier is not defined and the corresponding cells are excluded from
the training data; at inference such cells do not require any special
treatment, since $\delta_K^{\mathrm{NN}}=0$ for any predicted multiplier.}

The model is trained using mini-batch stochastic optimization with the
\texttt{AdamW} optimizer. The training is started with a learning rate of
$10^{-3}$ and a small weight decay of $10^{-5}$. To escape shallow local
minima, a cosine annealing schedule
\CHANGEF{with warm restarts is used: within each cycle the
learning rate decreases from its initial value following a cosine profile and is
then reset to that value at the beginning of the next cycle. The first cycle
lasts $T_0=50$ epochs and each subsequent cycle is twice as long as the previous
one ($T_{\mathrm{mult}}=2$), so that restarts occur after $50$, $150$, $350$,
$\dots$ epochs. To reduce the risk of exploding gradients and to enhance the
stability of the training, the Euclidean norm of the gradient is clipped so that
it does not exceed $1.0$. Additionally, early stopping is used to avoid
unnecessary computations and to prevent overfitting: if the validation loss does
not decrease by at least $10^{-5}$ over $50$ consecutive epochs, the training is
terminated.}
An overview of the hyperparameters used is
\CHANGEF{given}
in Tab.~\ref{tab:hyperparams}.

\begin{table}[t]
\centering
\begin{tabular}{|l|c|}
\hline
Hyperparameter & Value \\
\hline
Learning Rate & $10^{-3}$ \\
Weight Decay & $10^{-5}$ \\
Optimizer & AdamW \\
Scheduler Type & Cosine Annealing \\
Initial Restart Period ($T_0$) & 50 \\
Restart Multiplier ($T_{\text{mult}}$) & 2 \\
Dropout Base Rate & 0.2 \\
Gradient Clipping &  1.0 \\
Early Stopping Patience & 50 epochs \\
\hline
\end{tabular}
\caption{Hyperparameter configuration.}
\label{tab:hyperparams}
\end{table}

After the SUPG solution has been computed using the predicted stabilization
parameter, we apply a final post-processing step on cells $K \in \calG_h$
intersecting the outflow Dirichlet boundary. The Tabata solution $\uTAB$,
which satisfies the discrete maximum principle, is used as a reference: on
any cell $K \in \calG_h$ where $\uSUPG$ exhibits an overshoot relative to
$\uTAB$, the stabilization parameter $\deltaNN$ is multiplied by a factor
of $10$ and the SUPG system is resolved. The procedure is repeated until
no overshoots remain or until a maximum of five iterations is reached.

\section{Numerical studies for the Hemker problem}
\label{sec:hemker}

The Hemker problem is a standard benchmark problem for steady-state
convection-diffusion equations introduced in \cite{Hem96} and here we use the
following setting.

\begin{example}\label{ex:hemker}
Problem \eqref{steady-strong} is considered with 
$$
   \Omega = \big((-3,9)\times(-3,3)\big) \setminus 
   \{(x,y)\in{\mathbb R}^2\,;\,\,x^2+y^2 \le 1 \}\,,
$$
$\Gamma^N=((-3,9)\times\{-3\})\cup(\{9\}\times[-3,3])\cup((-3,9)\times\{3\})$,
$\Gamma^D=\partial\Omega\setminus\overline{\Gamma^N}$, $\varepsilon=10^{-4}$,
$\bb=(1,0)^\top$, $f=0$, $g=0$, and
\begin{align*}
   u_b(-3,y) &= 0\qquad\forall\,\,y\in(-3,3)\,,\\
   u_b(x,y) &= 1\qquad\forall\,\,(x,y)\in{\mathbb R}^2\,\,\,
   \mbox{\rm with}\,\,\,x^2+y^2=1\,.
\end{align*}
\end{example}

\begin{figure}[t!]
\begin{center}
\includegraphics[height=0.28\textheight]{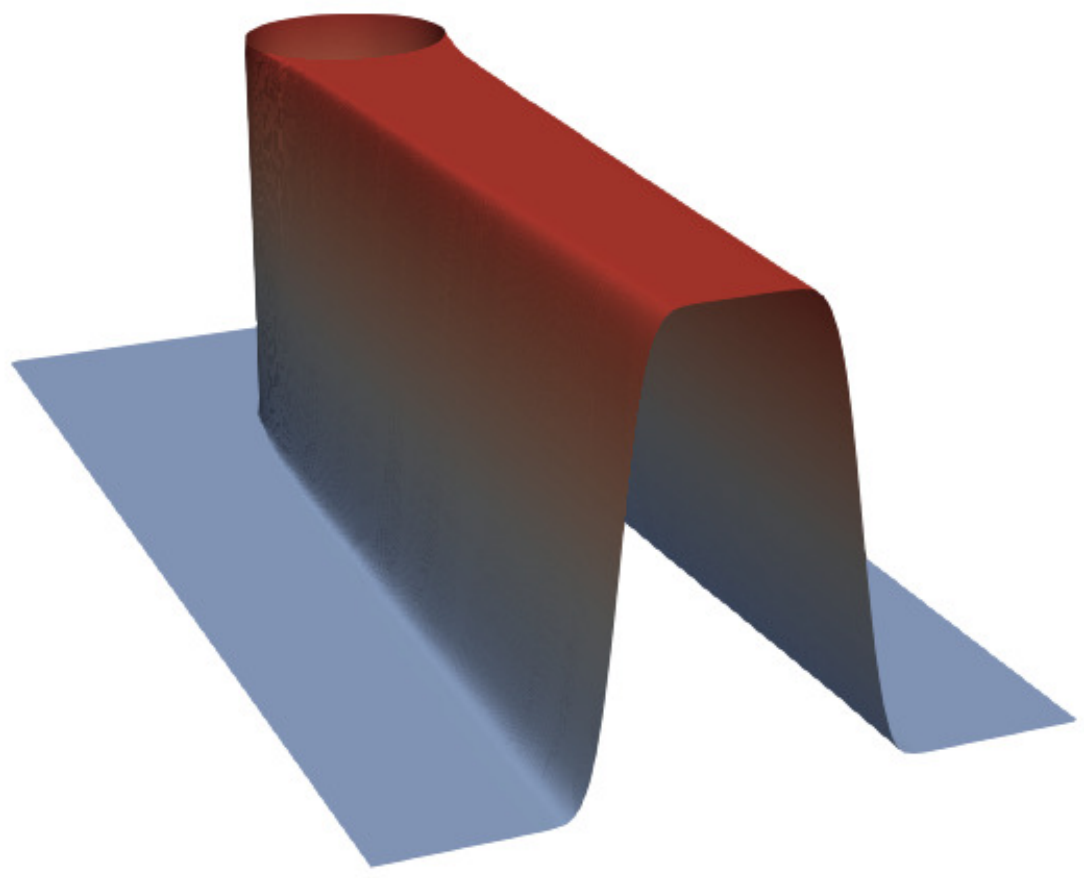}\\[2mm]
\includegraphics[height=0.28\textheight]{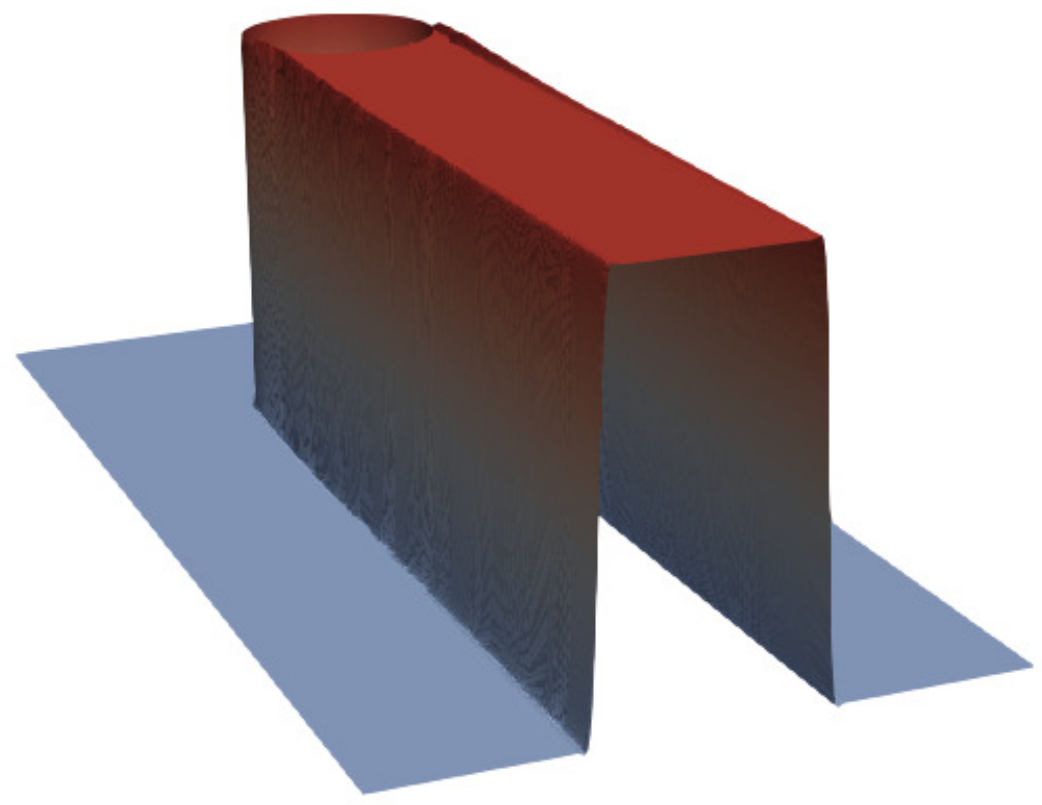}\\[2mm]
\includegraphics[height=0.28\textheight]{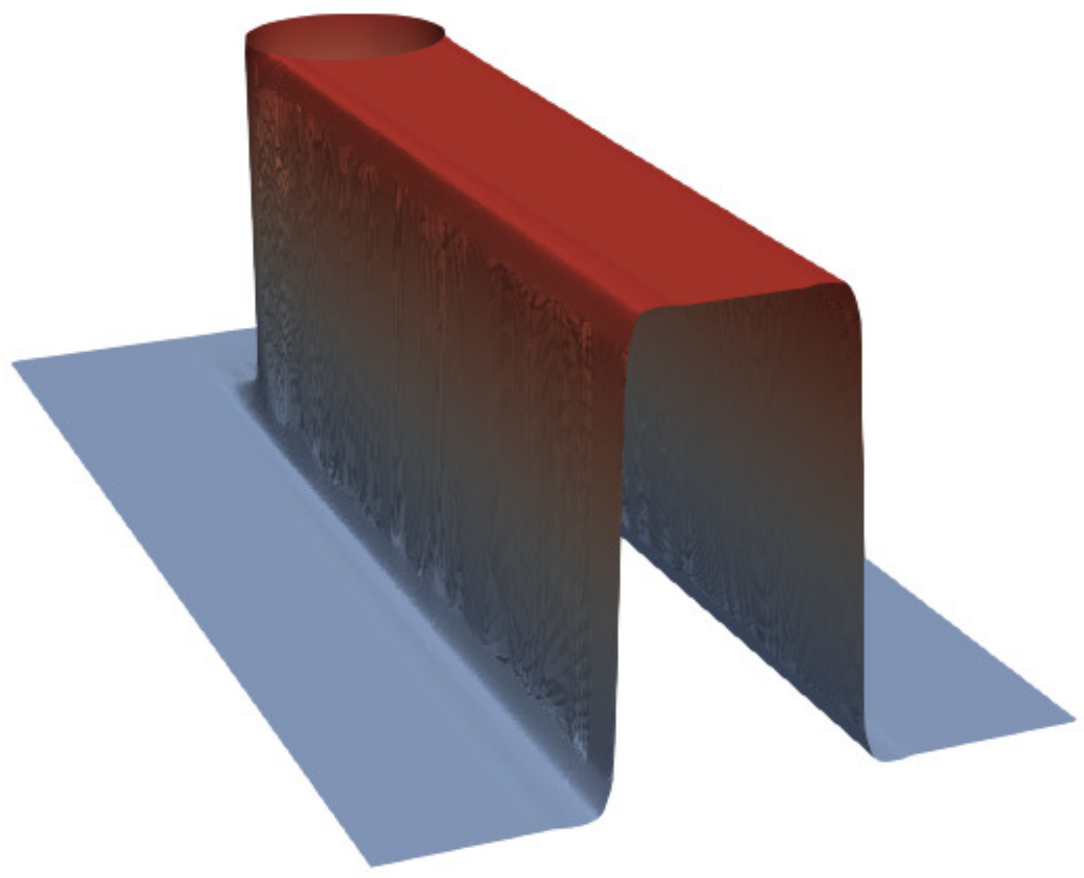}
\end{center}
\caption{Example~\ref{ex:hemker}: solutions on level~5 obtained using the
Tabata method (top), the SUPG method with $\deltaSTD$ and $\deltaBND$ (middle),
and the SUPG method with $\deltaNN$ (bottom).}\label{fig:hemker_sol}
\end{figure}

The Dirichlet boundary condition $u=0$ prescribed at $x=-3$ is transported by
the convection $\bb$ to the boundary part formed by the unit circle, where the
Dirichlet boundary condition $u=1$ is considered. This creates an exponential
boundary layer at the left-hand side of the circle, which continues to the top
and bottom of the circle. From there, two interior layers start in the direction
of the flow field, see the numerical results in Fig.~\ref{fig:hemker_sol}.
Note that the solution $u$ of Example~\ref{ex:hemker} attains values from the
interval $[0,1]$ only, which is a consequence of the maximum principle, see
\cite[Theorem 2.21]{BJK25}.

The setting of Example~\ref{ex:hemker} was considered in \cite{ACF+11} to
assess various numerical methods for problem \eqref{steady-strong}. To this
end, a reference solution was employed which was computed using the Galerkin
finite element method with the $Q_1$ finite element and a fine grid containing
$48 252 416$ degrees of freedom. This enabled to compare the accuracy of
various methods along various cuts of the solution. Here, we will concentrate
on the approximation of the interior layer at $y=1$ along $x=4$. In
\cite{ACF+11}, the width of the layer was defined as the length of the interval
$[y_0,y_1]$ in which the approximate solution $u_h$ falls from 0.9 to 0.1,
i.e.,
\begin{equation}\label{eq:y0_y1}
   u_h(4,y) > 0.9\quad\forall\,\,y\in(0,y_0)\,,\qquad
   u_h(4,y) < 0.1\quad\forall\,\,y\in(y_1,3)\,.
\end{equation}
For the reference solution, the layer width is 0.0723. We will compute the
absolute value of the difference between the layer width $y_1-y_0$ and this
reference value to compare various methods.
\CHANGEF{The values $y_0$ and $y_1$ in~\eqref{eq:y0_y1} are determined from the
values of $u_h$ at $100000$ equidistant points along $\{4\}\times[0,3]$. If the
thresholds are crossed several times (which never happened in the computations
presented in this paper), then the first crossing of $0.9$ and the last
crossing of $0.1$ are selected.}
Further quantities of interest
considered in this paper are the magnitudes of the undershoots and overshoots
which will be computed both globally and along the line $x=4$. They are defined
by
\begin{alignat}{4}
   &\min{}_{\mathrm{g}}\,&&=-\min_{\overline\Omega}\,u_h\,,\qquad
   &&\max{}_{\mathrm{g}}\,&&=\max_{\overline\Omega}\,u_h-1\,,
   \label{eq:min_max_global}\\
   &\min{}_{\mathrm{l}}&&=-\min_{[-3,3]}\,u_h(4,\cdot)\,,\qquad
   &&\max{}_{\mathrm{l}}&&=\max_{[-3,3]}\,u_h(4,\cdot)-1\,.\label{eq:min_max_local}
\end{alignat}
Then $\min{}_{\mathrm{g}}\ge0$, $\max{}_{\mathrm{g}}\ge0$,
$\min{}_{\mathrm{g}}\ge\min{}_{\mathrm{l}}$, and
$\max{}_{\mathrm{g}}\ge\max{}_{\mathrm{l}}$. If $\min{}_{\mathrm{l}}\le0$, then
there is no undershoot at $x=4$ but $u_h$ is positive along this cutline.
Similarly, if $\max{}_{\mathrm{l}}<0$, then there is no overshoot at $x=4$ but
$u_h(4,\cdot)<1$.

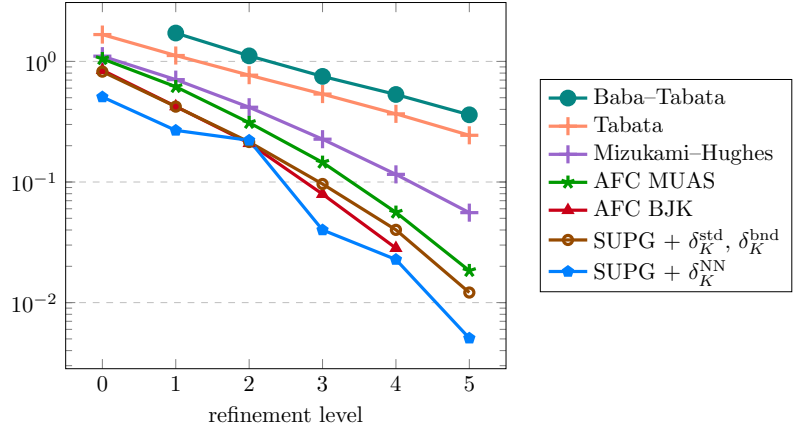
\begin{figure}[t!]
\begin{center}
\begin{tikzpicture}[scale=0.85]
\begin{axis}[
legend pos=south east, xlabel = refinement level, 
legend cell align ={left},
legend style={nodes={scale=1, transform shape},at={(1.65,0.5)},anchor=east},
xtick={0,1,2,3,4,5},
ymode=log,
ymajorgrids=true,
grid style=dashed]
\addplot[color=sea_green, line width = 0.5mm, solid, mark=oplus*, mark options={solid}, mark size=3pt]
coordinates{
(1,1.71879)
(2,1.11171)
(3,7.512599999999999e-01)
(4,5.3181e-01)
(5,3.6084e-01)
};
\addlegendentry{Baba--Tabata}
\addplot[color=salmon, line width = 0.5mm, solid, mark=+, mark options={solid}, mark size=4pt]
coordinates{
(0,1.669350e+00)
(1,1.119000e+00)
(2,7.688700e-01)
(3,5.360100e-01)
(4,3.674400e-01)
(5,2.439300e-01)
};
\addlegendentry{Tabata}
\addplot[color=amethyst, line width = 0.5mm, solid, mark=+, mark options={solid}, mark size=4pt]
coordinates{
(0,1.105260e+00)
(1,7.061400e-01)
(2,4.167900e-01)
(3,2.264400e-01)
(4,1.157100e-01)
(5,5.571000e-02)
};
\addlegendentry{Mizukami--Hughes}
\addplot[color=dark_green, line width = 0.5mm, solid, mark=star, mark options={solid}, mark size=3pt]
coordinates{
(0,1.05039)
(1,0.61635)
(2,0.3105)
(3,0.14541)
(4,0.05601)
(5,0.018390)
};
\addlegendentry{AFC MUAS}
\addplot[color=crimson, line width = 0.5mm, solid, mark=triangle*, mark options={solid}, mark size=2pt]
coordinates{
(0,0.84585)
(1,0.42339)
(2,0.20979)
(3,0.07938)
(4,0.02829)
};
\addlegendentry{AFC BJK}
\addplot[color=brown, line width = 0.5mm, solid, mark=o, mark options={solid}, mark size=2pt]
coordinates{
(0,8.226900e-01)
(1,4.217700e-01)
(2,2.146800e-01)
(3,9.624000e-02)
(4,4.011000e-02)
(5,1.212000e-02)
};
\addlegendentry{SUPG + $\deltaSTD$, $\deltaBND$}
\addplot[color=azure, line width = 0.5mm, solid, mark=pentagon*, mark options={solid}, mark size=2pt]
coordinates{
(0,0.50694)
(1,0.26799)
(2,0.2217)
(3,0.04002)
(4,0.02277)
(5,0.00507)
};
\addlegendentry{SUPG + $\deltaNN$}
\end{axis}
\end{tikzpicture}
\end{center}
\caption{Example~\ref{ex:hemker}: differences between the layer widths of the computed solutions and the reference layer width. 
}\label{fig:hemker_layer}
\end{figure}

\begin{figure}[t!]
\centerline{\includegraphics[width=0.7\textwidth]{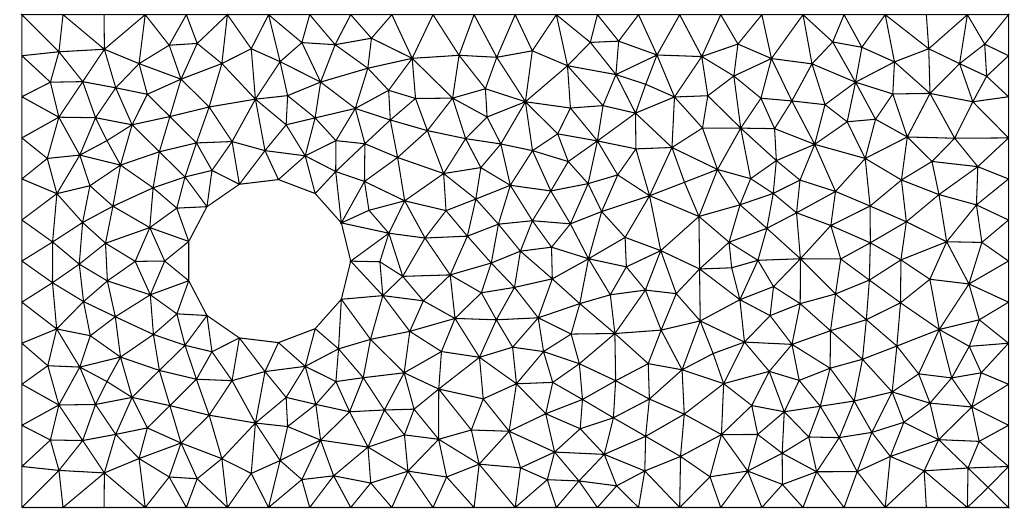}}
\caption{Example~\ref{ex:hemker}: initial triangulation (level 0).}
\label{fig:hemker_grid}
\end{figure}

Recently, the above defined layer width of Example~\ref{ex:hemker} was used in
\cite{JKP23} to assess several methods satisfying the discrete maximum
principles. Some of these results are presented in Fig.~\ref{fig:hemker_layer}
and compared with the approaches considered in this paper. We used the same
sequence of triangulations as in \cite{JKP23}. The initial triangulation (level~0) is depicted in Figure~\ref{fig:hemker_grid}. It was generated with 
\textsc{Gmsh} \cite{GR09}, is of Delaunay type, and consists of 743 mesh cells.
Finer triangulations (levels~1--5) were obtained by successive red refinement
and by adjusting the vertices near the circular part of the boundary.

The numerical results presented in \cite{JKP23} included linear upwinding
techniques leading to similar layer widths from which we present the results of
the Baba--Tabata method \cite{BT81}, the Mizukami--Hughes method \cite{MH85}
with the improvements \cite{Kno06}, which is a nonlinear method of upwind type
and is also presented in Fig.~\ref{fig:hemker_layer}, and a group of nonlinear
algebraic flux correction schemes which belong to most efficient approaches for
solving convection-dominated convection--diffusion problems according to
\cite{BJK23,BJK25}.
From this group, Fig.~\ref{fig:hemker_layer} shows the
results for the AFC scheme with the BJK limiter \cite{BJK17}, which gave the
smallest errors of the layer width in \cite{JKP23}, and for the MUAS method
\cite{JK21} where the errors were among the largest in the group of AFC
schemes. We observe that, among the methods from \cite{JKP23}, the largest 
smearing of the interior layer appears for the method of Baba and Tabata.
On level~0, the solution is even so smeared that $u_h(4,\cdot)<0.9$ so that the
layer width cannot be determined using \eqref{eq:y0_y1}. Therefore, the
corresponding result is missing in Fig.~\ref{fig:hemker_layer}. One may also
notice that Fig.~\ref{fig:hemker_layer} does not contain a result for the AFC
scheme with the BJK limiter on level~5, which is because the corresponding
nonlinear problem could not be solved. The difficulties with the convergence of
solvers represent an important bottleneck of most nonlinear methods.

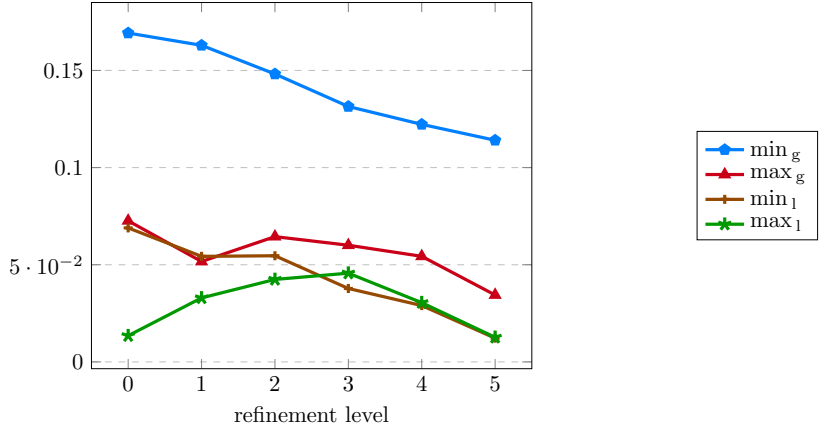
\begin{figure}[t!]
\begin{center}
\begin{tikzpicture}[scale=0.85]
\begin{axis}[
legend pos=south east, xlabel = refinement level, 
legend cell align ={left},
legend style={nodes={scale=1, transform shape},at={(1.65,0.5)},anchor=east},
xtick={0,1,2,3,4,5},
ymajorgrids=true,
grid style=dashed]

\addplot[color=azure, line width = 0.5mm, solid, mark=pentagon*, mark options={solid}, mark size=2pt]
coordinates{
(0,1.692537e-01)
(1,1.629293e-01)
(2,1.481734e-01)
(3,1.314471e-01)
(4,1.222484e-01)
(5,1.140726e-01)
};
\addlegendentry{$\min{}_{\mathrm{g}}$}

\addplot[color=crimson, line width = 0.5mm, solid, mark=triangle*, mark options={solid}, mark size=2pt]
coordinates{
(0,7.261948e-02)
(1,5.166819e-02)
(2,6.447137e-02)
(3,6.005154e-02)
(4,5.434293e-02)
(5,3.443813e-02)
};
\addlegendentry{$\max{}_{\mathrm{g}}$}

\addplot[color=brown, line width = 0.5mm, solid, mark=+, mark options={solid}, mark size=2pt]
coordinates{
(0,6.896397e-02)
(1,5.432348e-02)
(2,5.463392e-02)
(3,3.779404e-02)
(4,2.900196e-02)
(5,1.221454e-02)
};
\addlegendentry{$\min{}_{\mathrm{l}}$}

\addplot[color=dark_green, line width = 0.5mm, solid, mark=star, mark options={solid}, mark size=3pt]
coordinates{
(0,1.348955e-02)
(1,3.294879e-02)
(2,4.244250e-02)
(3,4.567257e-02)
(4,3.048845e-02)
(5,1.267153e-02)
};
\addlegendentry{$\max{}_{\mathrm{l}}$}

\end{axis}
\end{tikzpicture}
\end{center}
\caption{Example~\ref{ex:hemker}: measures for undershoots and overshoots
defined by \eqref{eq:min_max_global} and \eqref{eq:min_max_local} for solutions
of the SUPG method with $\deltaSTD$ and $\deltaBND$.
}\label{fig:min_max_SUPG_out}
\end{figure}

\begin{figure}[t!]
\begin{center}
\begin{tikzpicture}[scale=0.85]
\begin{axis}[
legend pos=south east, xlabel = refinement level, ylabel = {Error * $10^{-5}$},
legend cell align ={left},
legend style={nodes={scale=1, transform shape},at={(1.65,0.5)},anchor=east},
xtick={0,1,2,3,4,5},
ymajorgrids=true,
grid style=dashed]

\addplot[color=azure, line width = 0.5mm, solid, mark=pentagon*, mark options={solid}, mark size=2pt]
coordinates{
(0,4.9)
(1,3.5)
(2,3.8)
(3,3.5)
(4,4.9)
(5,5)
};
\addlegendentry{$\min{}_{\mathrm{g}}$}

\addplot[color=crimson, line width = 0.5mm, solid, mark=triangle*, mark options={solid}, mark size=2pt]
coordinates{
(0,3.1)
(1,0.1)
(2,5.4)
(3,3.9)
(4,4.3)
(5,4.6)
};
\addlegendentry{$\max{}_{\mathrm{g}}$}

\addplot[color=brown, line width = 0.5mm, solid, mark=+, mark options={solid}, mark size=2pt]
coordinates{
(0,4.9)
(1,3.4)
(2,0.2)
(3,2.1)
(4,4)
(5,2.7)
};
\addlegendentry{$\min{}_{\mathrm{l}}$}

\addplot[color=dark_green, line width = 0.5mm, solid, mark=star, mark options={solid}, mark size=3pt]
coordinates{
(0,1.6)
(1,3.7)
(2,1.9)
(3,3.7)
(4,2.1)
(5,4.6)
};
\addlegendentry{$\max{}_{\mathrm{l}}$}

\end{axis}
\end{tikzpicture}
\end{center}
\caption{Example~\ref{ex:hemker}: measures for undershoots and overshoots
defined by \eqref{eq:min_max_global} and \eqref{eq:min_max_local} for solutions
of the SUPG method with stabilization parameters $\deltaNN$ computed using the
neural network.
}\label{fig:min_max_SUPG_NN}
\end{figure}
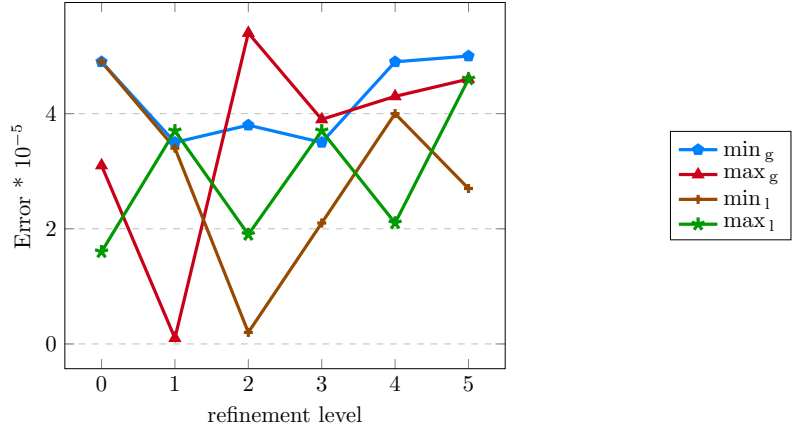

In addition to the methods already considered in \cite{JKP23},
Fig.~\ref{fig:hemker_layer} also shows the results for the Tabata method
mentioned in the preceding section which leads to a small improvement compared
to the Baba--Tabata method. On the other hand, we also observe from
Fig.~\ref{fig:hemker_layer} that the SUPG method with the stabilization
parameter $\deltaSTD$ on $K\in\calT_h\setminus\calG_h$ and $\deltaBND$ on
$K\in\calG_h$, cf.~\eqref{eq:delta_std} and \eqref{eq:delta0_bnd},
respectively, leads to a very sharp approximation of the interior layer,
comparable with the best method from \cite{JKP23}. However, in contrast to the
methods from \cite{JKP23} which satisfy the discrete maximum principle, the
SUPG method with $\deltaSTD$ and $\deltaBND$ leads to significant undershoots
and overshoots as Fig.~\ref{fig:min_max_SUPG_out} shows.
%
%
This deficiency can be significantly diminished by computing the stabilization
parameters using a neural network as described in the preceding section. Then
the undershoots and overshoots are
\CHANGEF{reduced by about three to four orders of magnitude: the quantities
defined in \eqref{eq:min_max_global} and \eqref{eq:min_max_local} do not exceed
$5\cdot 10^{-5}$ on any of the six levels, see Fig.~\ref{fig:min_max_SUPG_NN},
whereas they attain values up to $0.17$ for the standard parameters, see
Fig.~\ref{fig:min_max_SUPG_out}. Note the different scaling of the vertical
axes in the two figures. At the same time, the approximation of the interior
layer does not deteriorate: as Fig.~\ref{fig:hemker_layer} shows, the layer
width computed with $\deltaNN$ is at least as accurate as for the best
nonlinear methods or better than them.}
The solutions obtained using the three methods just discussed are depicted in
Fig.~\ref{fig:hemker_sol}. The Tabata method provides a solution without
undershoots and overshoots but with significantly smeared interior layers. The
SUPG method leads to a sharp approximation of the layers in both cases but if
the stabilization parameters are defined by $\deltaSTD$ and $\deltaBND$, clearly
visible undershoots and overshoots appear. On the other hand, if the
stabilization parameters are computed by the neural network, undershoots and
overshoots are essentially suppressed.

\CHANGEF{
Finally, we compare the computational cost of the approaches considered above.
Table~\ref{tab:times-iso} reports wall-clock times for the six refinement levels
of Example~1. The column SUPG contains the time for a single solve of~\eqref{eq:supg}
with the standard parameters $\delta_K^{\mathrm{std}}$, $\delta_K^{\mathrm{bnd}}$.
The column OPT contains the time needed by the a~posteriori parameter optimization
of~\cite{JKS11}. The last two columns contain the time of the complete pipeline
proposed here, i.e., the SUPG solve, the Tabata solve, the extraction of the local
features, the network inference and the SUPG solve with the predicted
stabilization parameters, without and with the post-processing step of
Section~3; the number in parentheses is the number of additional solves performed
by the post-processing. All computations reported in this section were performed on a laptop with a 13th generation Intel Core i7-13620H processor, using a single core and without any GPU
acceleration, so that the network inference was also carried out on the CPU.

\begin{table}[b]
  \centering
  \caption{Example 1: computational times on the isotropic triangulations. The
    number in parentheses is the number of additional SUPG solves performed in
    the post-processing step.}
  \label{tab:times-iso}
  \begin{tabular}{|c|r|r|r|r|}
    \hline
    \textbf{Level} & \textbf{SUPG} & \textbf{OPT}
      & $\boldsymbol{\delta_K^{\mathrm{NN}}}$ & $\boldsymbol{\delta_K^{\mathrm{NN}}}$\\
     & & & \textbf{(no post-proc.)} & \textbf{(post-proc.)}\\ \hline
     0 & 4.9054 s & 15.38 s & 16.937 s & 24.739 s (2)\\ \hline
    1 & 4.9364 s & 287.25 s & 17.093 s & 25.541 s (2)\\ \hline
    2 & 5.1931 s & 5392.34 s & 17.429 s & 27.878 s (2)\\ \hline
    3 & 5.9927 s & 40690.38 s & 19.282 s & 31.706 s (2)\\ \hline
    4 & 7.8749 s & 210803.18 s & 26.592 s & 42.775 s (2)\\ \hline
    5 & 20.2598 s & --- & 76.103 s & 102.782 s (1)\\ \hline
  \end{tabular}
\end{table}

The cost of the parameter optimization grows roughly by an order of magnitude per
refinement level, from $15$~s on level~0 to more than $2\cdot 10^{5}$~s on level~4,
and on level~5 the computation was no longer feasible. In contrast, the cost of the
proposed approach grows at essentially the same rate as a single SUPG solve: the
complete pipeline requires around $3$ times the time of one SUPG solve on all six
levels without post-processing. Consequently, the speedup with respect to the optimization increases from a
factor of~$17$ on level~1 to about $10^{3}$ on level~4. We emphasize that on the
coarsest level the optimization is still slightly cheaper than the network-based
approach; the advantage of the latter appears once the discrete
problems become large.

The difference between the last two columns is the cost of the post-processing
step, which consists of the additional SUPG solves listed in parentheses and agrees
with these solve times to within the measurement accuracy. On levels~0--4 two
additional solves were performed and on level~5 a single one, so that the limit of
five iterations was never reached. The post-processing is applied only on the cells
intersecting $\Gamma^{\mathrm{out}}$.

Since the training of the networks is performed offline and only once, its cost is
not included in Table~\ref{tab:times-iso}; the same trained models are used on all refinement levels
considered here.
}

To obtain even a better approximation of the interior layers, meshes with
an\-iso\-trop\-ic cells aligned with the convection are of advantage in the
layer regions. Unfortunately, it turns out that the approach presented in
Section~\ref{sec:ml_iso} does not work in a satisfactory way when applied to
anisotropic meshes. Therefore, some extensions and modifications of this
approach will be introduced in the next section.

\section{Computation of $\delta_K$ using machine learning on anisotropic meshes}
\label{sec:ml_aniso}
The strategy from isotropic case carries over to anisotropic meshes with only
minor modifications. The local quantities defined in \eqref{eq:u_quantities}
still allow us to detect layer regions and spurious oscillations. Similarly as
before, we treat interior cells $K\in\calT_h\setminus\calG_h$ and cells
$K\in\calG_h$ that intersect the outflow  boundary separately by training a
separate model for each region. We confine ourselves to scaled problems whose
solutions satisfy $\min_{\overline\Omega}u=0$ and 
$\max_{\overline\Omega}u\sim1$. This simplifies the definition of the input
values for the neural network.

For interior cells, we keep the MLP architecture from Section~\ref{sec:ml_iso},
but the three feature inputs \eqref{eq:delta_int} turn out to be insufficient
on anisotropic meshes. We therefore extend the input to the following nine 
features:
\begin{align*}
   &\Pe_K\,,\quad h_K^\parallel\,,\quad h_K^\perp\,,\quad 
   \uSUPG({\mathbf x}_K^1)\,,\quad \uSUPG({\mathbf x}_K^2)\,,\quad 
   \uSUPG({\mathbf x}_K^3)\,,\\
   &\|\nabla\uSUPG\|_{L^\infty(K)}^{}\,,\quad
   \|\uSUPG-\uTAB\|_{L^\infty(\Delta_K)}^{}\,,\quad
   \|\nabla(\uSUPG-\uTAB)\|_{L^\infty(\Delta_K)}^{}\,,
\end{align*}
where ${\mathbf x}_K^1$, ${\mathbf x}_K^2$, and ${\mathbf x}_K^3$ are the
vertices of
$K$,
\CHANGEF{numbered in such a way that 
$\uSUPG({\mathbf x}_K^1)\ge\uSUPG({\mathbf x}_K^2)\ge\uSUPG({\mathbf x}_K^3)$,
which prevents a dependence of the input data on the numbering of the
vertices in the triangulation.}
Apart from the input dimension, the loss functional and the
training protocol are identical to those described in
Section~\ref{sec:ml_iso}.

The boundary cells turned out to be more complicated. The same architecture,
even with an enlarged feature set, did not lead to a sufficiently accurate
approximation of $\deltaOPT$. Two factors are likely reasonable: (i) the
boundary behaviour of $\deltaOPT$ depends on more local features than for
interior cells and (ii) the relevant interaction between these features is not
known a priori. We therefore use twelve input features: in the above nine
features the P\'eclet number $\Pe_K$ is replaced by $\angle\bb_K$ introduced in
Section~\ref{sec:ml_iso} and we add the following three interactive features
$$
   \frac{h_K^\perp}{h_K^\parallel}\,,\qquad 
   h_K^\parallel\,\|\nabla\uSUPG\|_{L^\infty(K)}^{}\,,\qquad 
   \frac{\deltaSTD}{\deltaBND}\,.
$$
The MLP architecture for the model will be replaced by a tabular Transformer
inspired by the FT-Transformer \cite{gorishniy2021revisiting}. We do not
pre-specify which feature combinations are relevant. The self-attention
mechanism learns these interactions during training.

The model adapts a Transformer encoder to numerical tabular data. Given an
input vector $x\in\mathbb{R}^{12}$, a linear map first lifts $x$ to
$\mathbb{R}^{12d}$, 
\CHANGEF{where $d=64$ is the embedding dimension,} 
and reshapes the result into twelve vectors $t_1, \ldots, t_{12} \in
\mathbb{R}^d$, one per feature. A trainable summary vector \texttt{[CLS]} $t_0
\in \mathbb{R}^d$ is prepended to this collection and a trainable positional
shift is added to each $t_i$ for $i = 0, \ldots, 12$. The resulting thirteen
vectors form the input to the encoder. The role of $t_0$ is to act as a
designated slot from which the final prediction will be read; during encoding
it accumulates information from the feature vectors.

The encoder consists of $L=2$ identical blocks. Each block performs two
operations. The first is multi-head self-attention with $H=4$ heads. For each
head, three matrices  $W_Q, W_K, W_V \in
\mathbb{R}^{d \times d_h}$ with $d_h = d/H$ are introduced. These are the ordinary parameter matrices; their entries are initialized at random. They are adjusted by gradient descent during training in the same way as the weights of an MLP. From each input vector $t_i$ the head computes three derived vectors, 
\begin{equation*}
   q_i = W_Q^\top t_i, \qquad k_i = W_K^\top t_i, \qquad v_i = W_V^\top t_i,
\end{equation*}
called the query, key and value of $t_i$. Attention weights are then defined for every pair $(i,j)$ by
\begin{equation*}
   \alpha_{ij} = \frac{\exp\!\big(\langle q_i, k_j \rangle / \sqrt{d_h}\big)}
                       {\sum_{\ell=0}^{12}
                        \exp\!\big(\langle q_i, k_\ell \rangle / \sqrt{d_h}\big)},
   \qquad \sum_{j=0}^{12} \alpha_{ij} = 1,
\end{equation*}
and the head replaces each $t_i$ by the convex combination $\sum_{j=0}^{12} \alpha_{ij}\, v_j$. The matrices $W_Q$ and $W_K$ together determine which pairs of input vectors are considered similar, while $W_V$ determines what information is exchanged once the similarities have been resolved. The four heads carry out this construction in parallel with independent matrices, such that four different similarity measures act on the input simultaneously and their outputs are concatenated. Self-attention is the only operation in the model that lets the thirteen vectors exchange information.

The second operation  in each block is a feed-forward network of shape $d \to
4d \to d$ with a \texttt{GELU} activation, applied independently to each of the
thirteen vectors. Layer normalization and residual connections are applied
around both operations to stabilize training.

\begin{figure}[htbp]
\centering
\resizebox{\textwidth}{!}{%
\begin{tikzpicture}[%
    font=\sffamily\small,
    >={Stealth[length=2.2mm]},
    every node/.style={align=center},
    block/.style ={draw, rectangle, rounded corners=2pt, thick,
                   minimum width=2.4cm, minimum height=0.85cm, line width=0.5pt},
    tok/.style   ={draw, rectangle, rounded corners=2pt, thick,
                   minimum width=0.75cm, minimum height=0.65cm, line width=0.5pt,
                   font=\footnotesize},
    inp/.style    ={block, fill=blue!12},
    embed/.style  ={block, fill=teal!18},
    cls/.style    ={tok, fill=orange!35},
    feat/.style   ={tok, fill=teal!18},
    sublyr/.style ={draw, rectangle, rounded corners=2pt, thick,
                    minimum width=2.6cm, minimum height=0.9cm, line width=0.5pt,
                    fill=purple!12, font=\footnotesize},
    head/.style   ={block, fill=red!18, minimum width=2.4cm},
    outblk/.style ={block, fill=green!22, minimum width=1.6cm},
    flow/.style   ={->, thick, line width=0.5pt},
    encbox/.style ={draw=purple!60, dashed, rounded corners=4pt, thick,
                    inner sep=10pt, line width=0.5pt},
    lyr/.style    ={align=center, font=\footnotesize\itshape, text=purple!50!black}
]
\node[inp] (input) {Input features\\$x \in \mathbb{R}^{12}$};
\node[embed, right=0.8cm of input] (tokenizer)
    {Feature tokenizer\\\footnotesize Linear $12 \to 12 d$};
\draw[flow] (input) -- (tokenizer);

\coordinate (trow) at ($(tokenizer.east)+(1.7,0)$);
\node[cls]  (t0)  at (trow)            {[CLS]};
\node[feat] (t1)  at ($(t0)+(0.85,0)$) {$t_1$};
\node[feat] (t2)  at ($(t1)+(0.85,0)$) {$t_2$};
\node[feat] (t3)  at ($(t2)+(0.85,0)$) {$\cdots$};
\node[feat] (t12) at ($(t3)+(0.85,0)$) {$t_{12}$};
\node[font=\footnotesize\itshape, above=0.15cm of t2, xshift=0.4cm]
    {Token sequence (length 13)};
\draw[flow] (tokenizer.east) -- ($(t0.west)+(-0.05,0)$);

\coordinate (busC) at ($(t3.south)+(0,-0.6)$);
\coordinate (busL) at (t0.south |- busC);
\coordinate (busR) at (t12.south |- busC);
\foreach \i in {0,1,2,3,12} \draw[thin] (t\i.south) -- (t\i.south |- busC);
\draw[thick, line width=0.6pt] (busL) -- (busR);

\coordinate (encIn) at ($(busC)+(0,-0.6)$);
\draw[flow] (busC) -- (encIn);

\node[sublyr, below=0.0cm of encIn, anchor=north] (mha)
    {Multi-head self-attention\\\footnotesize ($H{=}4$)};
\node[sublyr, below=0.6cm of mha] (ffn)
    {Feed-forward\\\footnotesize $d \to 4d \to d$};

\draw[flow] (mha) -- (ffn);

\node[encbox, fit=(mha)(ffn),
      label={[lyr]left:{$\times\,L{=}2$}}] (encblock) {};

\coordinate (encOut) at ($(ffn.south)+(0,-0.55)$);
\draw[flow] (ffn.south) -- (encOut);

\node[font=\footnotesize\itshape, left=0.1cm of encOut, xshift=-0.2cm] {Take [CLS]};

\node[head, right=0.7cm of encOut, anchor=west] (head1)
    {Linear $d{\to}64$\\\footnotesize \texttt{GELU}, dropout};
\node[head, right=0.5cm of head1] (head2)
    {Linear $64{\to}1$\\\footnotesize Sigmoid};
\node[outblk, right=0.5cm of head2]  (mainout) {Output};

\draw[flow] (encOut) -- (head1.west);
\draw[flow] (head1) -- (head2);
\draw[flow] (head2) -- (mainout);

\end{tikzpicture}%
}
\caption{Tabular Transformer architecture for the boundary cells.}
\label{fig:tabular-transformer}
\end{figure}
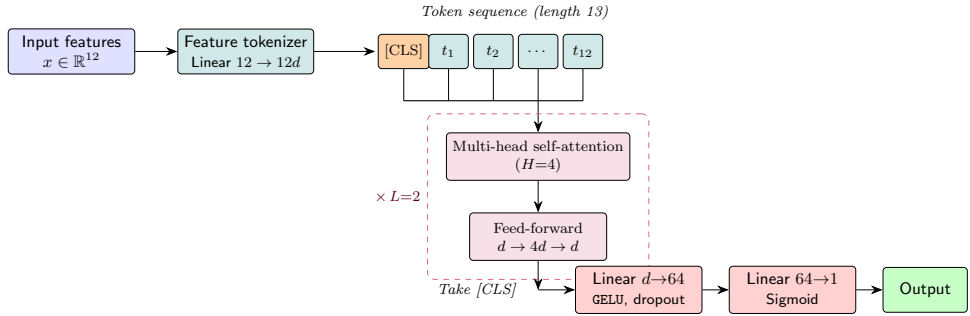

After the encoder, the twelve feature vectors are discarded and only the
summary vector is retained. It is passed through a final layer normalization
and a two-layer perceptron of shape $d \to 64 \to 1$ with \texttt{GELU}
activation and a sigmoid output.  The full architecture is shown in
Fig.~\ref{fig:tabular-transformer}.

The training protocol is similar to the one used for the isotropic model, but
several choices were adjusted to better fit the Transformer architecture. We
use the \texttt{AdamW} optimizer with a learning rate of $10^{-3}$ and a weight
decay of $10^{-5}$. The loss considered is the mean absolute error
(\texttt{MAE}) between the predicted and optimal stabilization parameter. The
learning rate is reduced by a factor of 0.5 whenever the validation loss has
not improved for 15 consecutive epochs (\texttt{ReduceLROnPlateau}). To prevent
the gradient from exploding, the gradient is clipped at $1$. Training runs for
$300$ epochs, with early stopping if the validation loss does not improve by at
least $10^{-5}$ for $30$ consecutive epochs. A dropout rate of $0.1$ is applied
inside the encoder layers. The complete hyperparameter configuration is listed
in  Tab.~\ref{tab:hyperparams_tt}.

\begin{table}[t]
\centering
\begin{tabular}{|l|c|}
\hline
Hyperparameter & Value \\
\hline
Embedding dimension $d$ & 64 \\
Number of encoder layers $L$ & 2 \\
Number of attention heads $H$ & 4 \\
Feed-forward width & $4d$ \\
Dropout & 0.1 \\
Loss & MAE \\
Optimizer & AdamW \\
Learning rate & $10^{-3}$ \\
Weight decay & $10^{-5}$ \\
LR scheduler & ReduceLROnPlateau \\
Scheduler factor & 0.5 \\
Scheduler patience & 15 epochs \\
Gradient clipping & 1.0 \\
Maximum epochs & 300 \\
Early stopping tolerance & $10^{-5}$ \\
Early stopping patience & 30 epochs \\
Ensemble size & 5 \\
\hline
\end{tabular}
\caption{Hyperparameter configuration for the tabular Transformer used on the
boundary cells.}
\label{tab:hyperparams_tt}
\end{table}

After the SUPG solution is computed, we apply the same post-processing step as
in the isotropic case on the outflow boundary cells which exhibit overshoots,
see the end of Section~\ref{sec:ml_iso}.

\section{Numerical studies for anisotropic meshes}
\label{sec:numerics_aniso}

In this section, we present numerical results computed on anisotropic meshes
obtained using the ANGENER software \cite{angener}. Starting from an isotropic
triangulation and the corresponding solution of the SUPG method with the
standard stabilization parameter \eqref{eq:delta_std}, ANGENER constructs an
anisotropic triangulation using the anisotropic mesh adaptation (AMA)
algorithm \cite{Dol98,DF04}. The resulting meshes have strongly stretched cells
along interior and boundary layers.

\begin{figure}[t!]
\centerline{\includegraphics[width=0.7\textwidth]{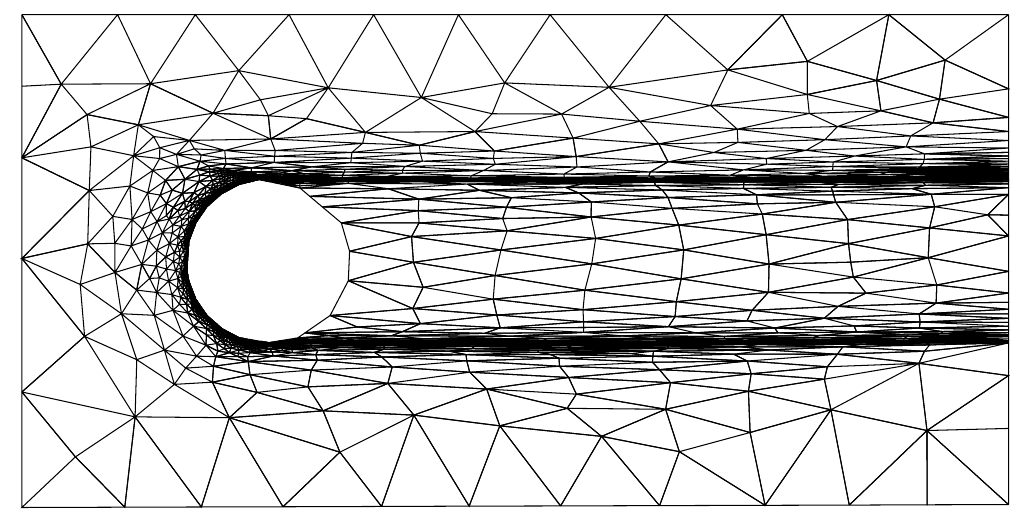}}
\caption{Example~\ref{ex:hemker}: triangulation obtained by anisotropic mesh adaptation.}
\label{fig:hemker_grid_aniso} 
\end{figure}

\begin{figure}[t!]
\centerline{\includegraphics[width=\textwidth]{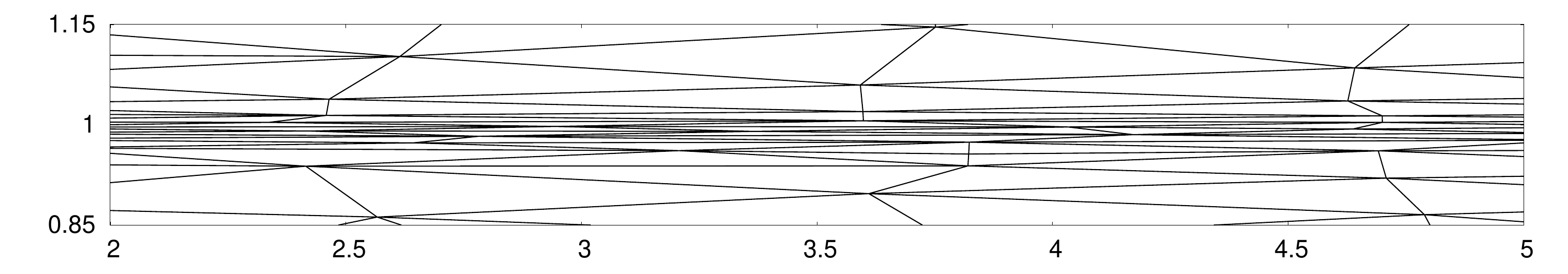}}
\caption{Example~\ref{ex:hemker}: zoom of the anisotropic mesh along the
interior layer at $y=1$. Notice the different scaling of the coordinate axes.}
\label{fig:hemker_grid_aniso_zoom} 
\end{figure}

\begin{figure}[t!]
\begin{center}
\includegraphics[height=0.35\textheight]{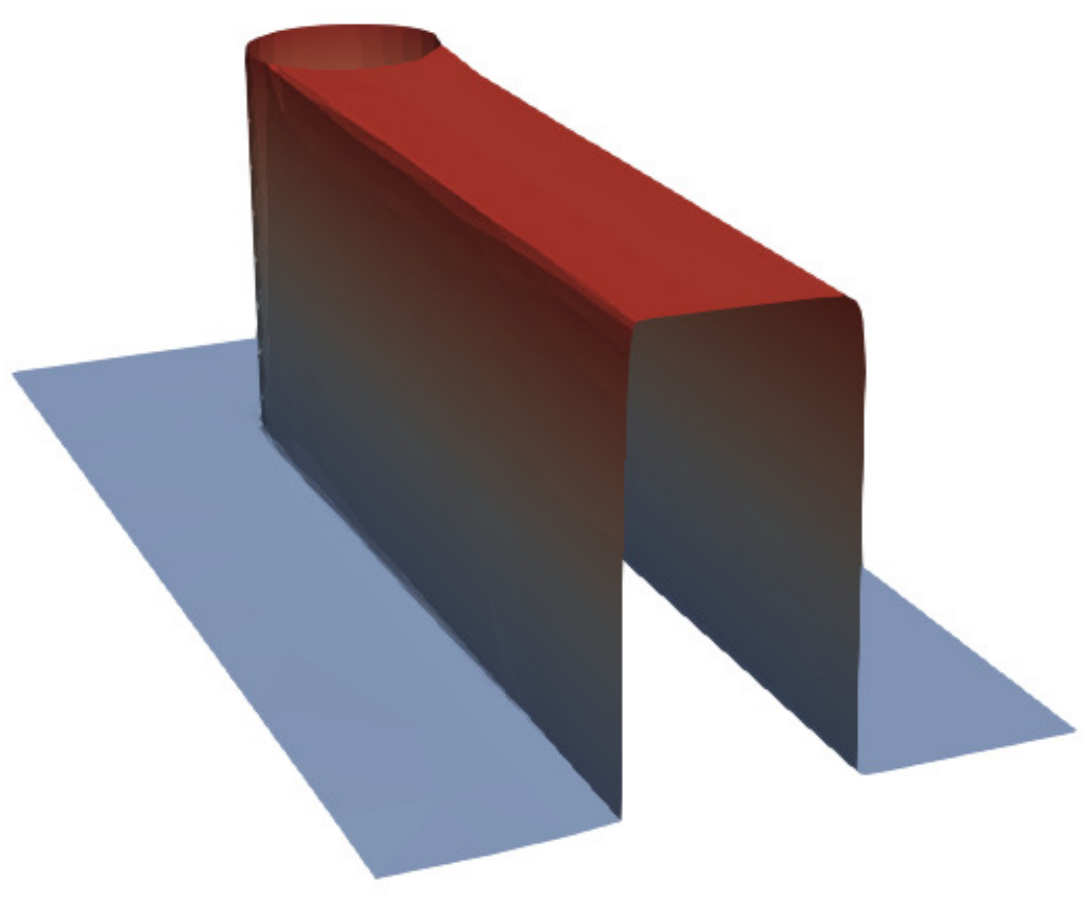}
\end{center}
\caption{Example~\ref{ex:hemker}: solution of the SUPG method with $\deltaNN$
obtained on the anisotropic mesh from Fig.~\ref{fig:hemker_grid_aniso}.}
\label{fig:hemker_sol_aniso}
\end{figure}

Let us start with the Hemker problem investigated in Section~\ref{sec:hemker}.
The triangulation obtained using ANGENER is shown in
Fig.~\ref{fig:hemker_grid_aniso} and a zoom of a part of the triangulation 
in the region of the interior layer at $y=1$ is depicted in
Fig.~\ref{fig:hemker_grid_aniso_zoom}. The solution of the SUPG method with
the stabilization parameters $\deltaNN$ computed using the approach described in
Section~\ref{sec:ml_aniso} approximates the reference layer width 0.0723 with
the error $0.0084$, which is a comparable accuracy as obtained using the SUPG
method with $\deltaNN$ on level~5 of the isotropic triangulations. However, 
this isotropic triangulation consists of 760832 mesh cells whereas the
anisotropic triangulation contains 5250 mesh cells only. Consequently, the
computation on the anisotropic triangulation is significantly cheaper.
Moreover, the quantities defined in \eqref{eq:min_max_global} and
\eqref{eq:min_max_local} take now the following values
$$
   \min{}_{\mathrm{g}}=0.003876\,,\quad
   \max{}_{\mathrm{g}}=0.003850\,,\quad
   \min{}_{\mathrm{l}}=0.002511\,,\quad
   \max{}_{\mathrm{l}}=0.003657\,,
$$
\CHANGEF{The solution is depicted in Fig.~\ref{fig:hemker_sol_aniso} and is}
more accurate than the solutions obtained on isotropic meshes
\CHANGEF{with a comparable number of degrees of freedom.}

Next let us consider another classical benchmark problem proposed
in~\cite{HMM86}, whose solution contains an interior layer and two exponential
boundary layers.

\begin{figure}[t!]
\centerline{\includegraphics[width=0.4\textwidth]{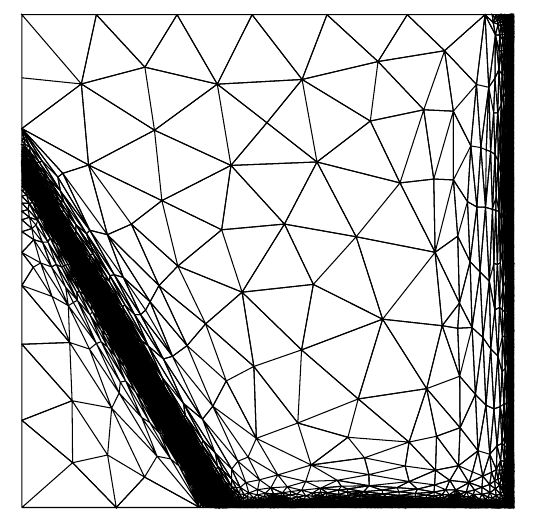}
\includegraphics[width=0.4\textwidth]{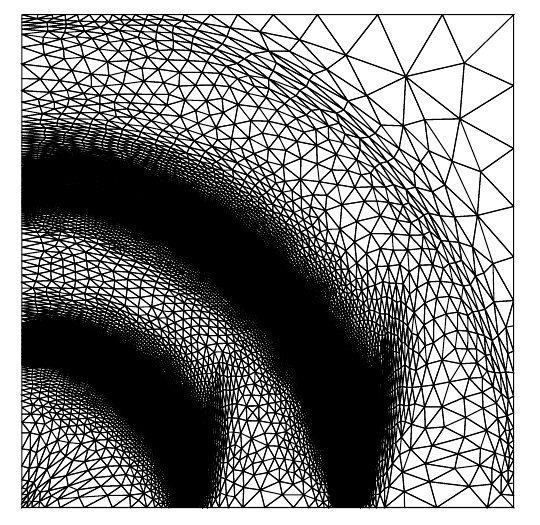}
}
\caption{Triangulations obtained by anisotropic mesh
adaptation used for approximating the solutions of Example~\ref{ex1} (left)
and  Example~\ref{ex2} (right).}
\label{fig:aniso_meshes} 
\end{figure}

\begin{figure}[t!]
\begin{center}
\includegraphics[width=0.45\textwidth]{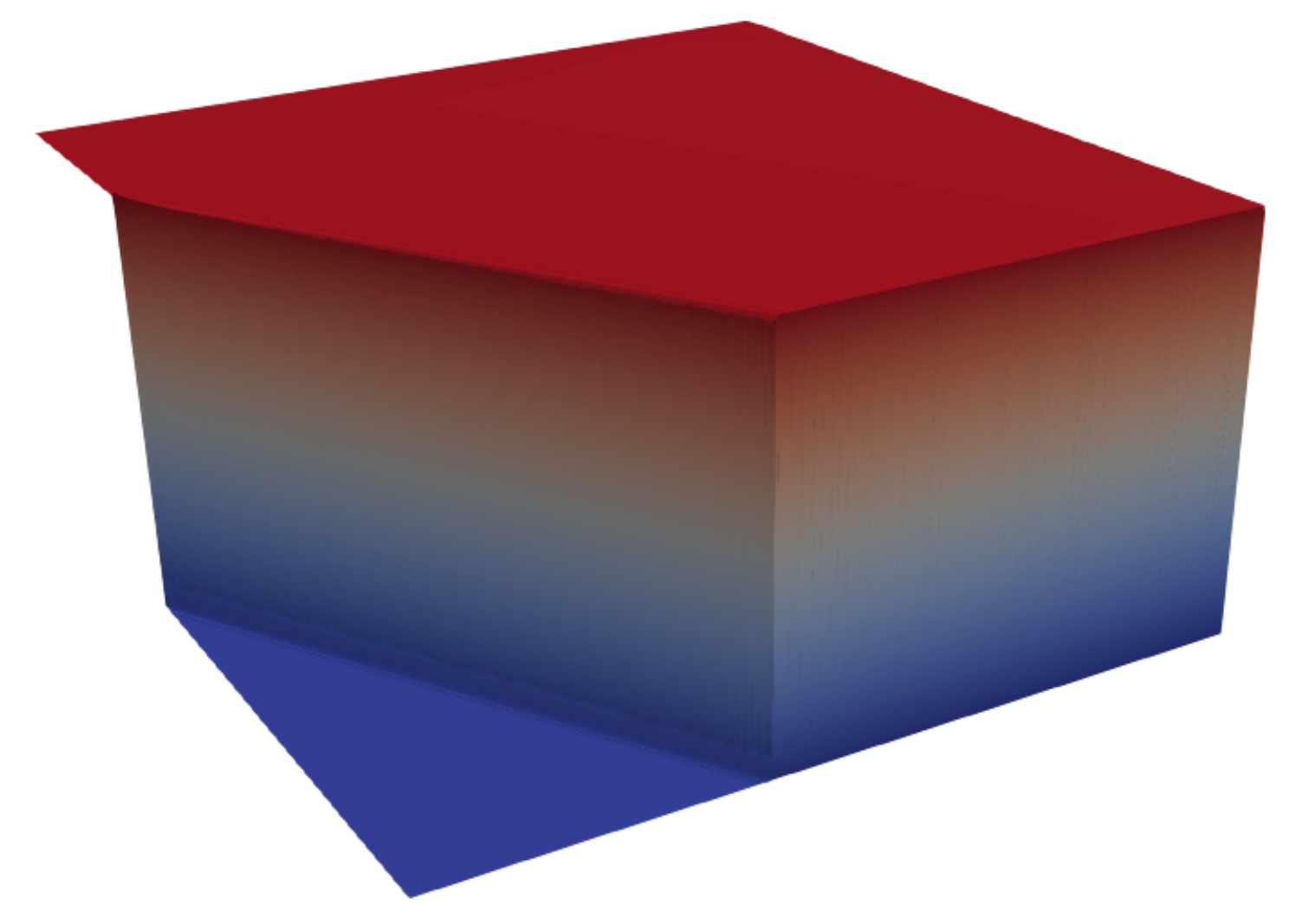}
\hspace*{3mm}
\includegraphics[width=0.45\textwidth]{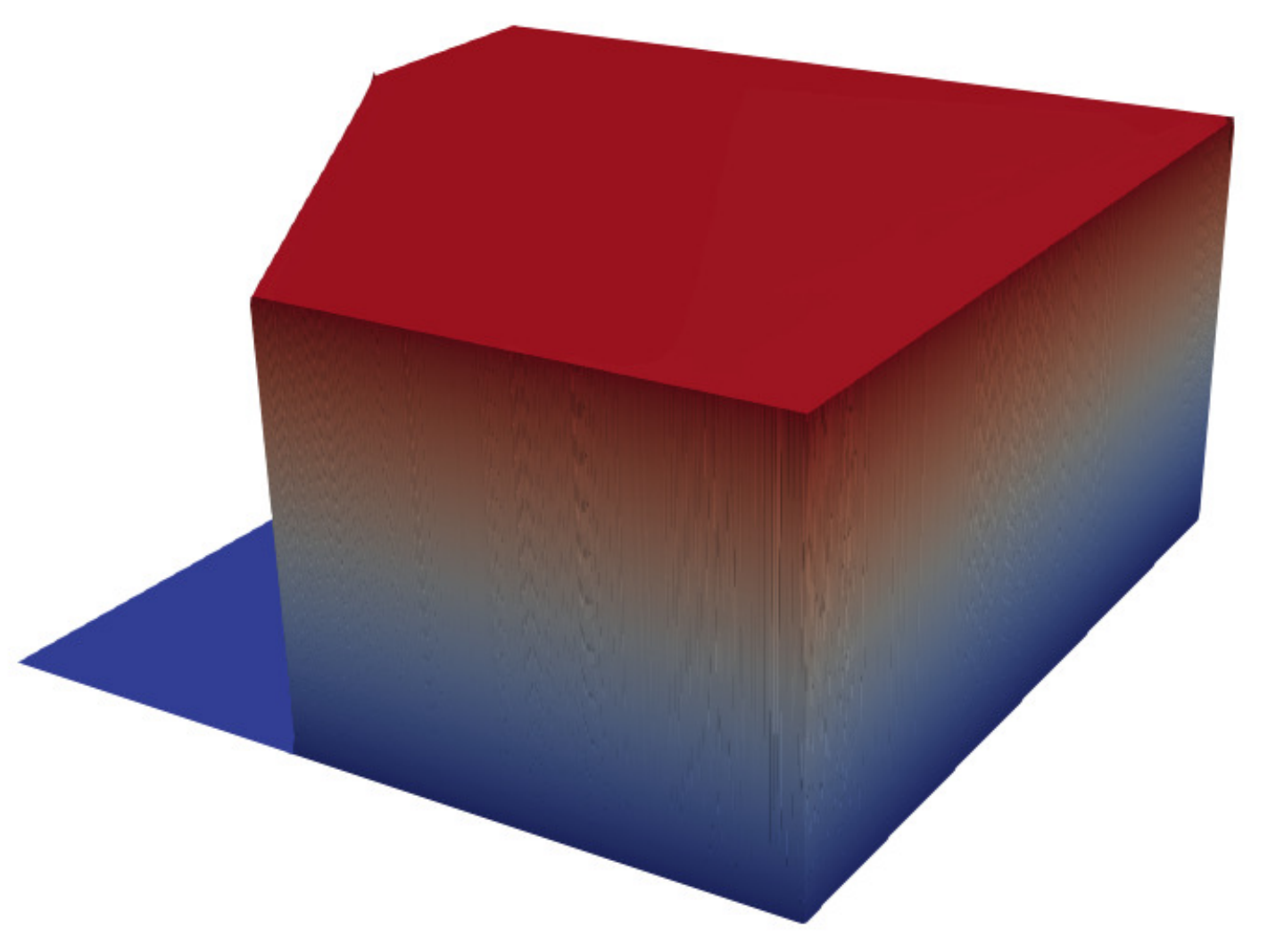}
\end{center}
\caption{Example~\ref{ex1}: two views of the solution of the SUPG method with
$\deltaNN$ obtained on the anisotropic mesh from Fig.~\ref{fig:aniso_meshes}
(left).}
\label{fig:ex1_sol}
\end{figure}

\begin{example}\label{ex1}
Problem \eqref{steady-strong} is considered with $\Omega = (0,1)^2$,
$\Gamma^D=\partial\Omega$, $\Gamma^N=\emptyset$, $\varepsilon=10^{-8}$,
$\bb=(\cos(-\pi/3),\sin(-\pi/3))^\top$, $f=0$, and
$$
u_b(x,y) = \left\{ \begin{array}{ll}
0 & \quad\mbox{for $x=1$ or $y\le0.7$,}\\
1 & \quad\mbox{else.}
\end{array}\right.
$$ 
\end{example}
The triangulation used for this example is shown in Fig.~\ref{fig:aniso_meshes}
(left). The quality of the corresponding SUPG solution with stabilization
parameters computed using the approach of Section~\ref{sec:ml_aniso} is
illustrated in Fig.~\ref{fig:ex1_sol}.
\CHANGEF{All layers are approximated without visible smearing and the
unphysical values are reduced by more than an order of magnitude in comparison
with the standard parameters, see Table \ref{tab:aniso-ex23}.}
\CHANGEF{
\begin{table}[htbp]
  \centering
  \caption{Examples \ref{ex1} and \ref{ex2}: measures for undershoots and overshoots defined
    by \ref{eq:min_max_global} for the solutions on the anisotropic meshes, computed with
    the standard parameters $\delta_K^{\mathrm{std}}$
    and with the predicted parameters $\delta_K^{\mathrm{NN}}$. The numbers in
    parentheses are the numbers of additional solves performed in the
    post-processing step.}
  \label{tab:aniso-ex23}
  \begin{tabular}{|c|c|r|r|}
    \hline
    \textbf{Problem} & \textbf{Quantity}
      & $\boldsymbol{\delta_K^{\mathrm{std}}}$
      & $\boldsymbol{\delta_K^{\mathrm{NN}}}$\\ \hline
    Example 2 & $\min_{\mathrm g}$ & 0.0438 & 0.0038 (2)\\ \hline
    Example 2 & $\max_{\mathrm g}$ & 0.3148 & 0.0143 (2)\\ \hline
    Example 3 & $\min_{\mathrm g}$ & 0.0010 & 0.0001 (0)\\ \hline
    Example 3 & $\max_{\mathrm g}$ & 0.1560 & 0 (0)\\ \hline
  \end{tabular}
\end{table}
}

The last test problem we will consider was proposed in \cite{KLS19} and is
defined by the following data.

\begin{example}\label{ex2}
Problem \eqref{steady-strong} is considered with $\Omega = (0,1)^2$,
$\Gamma^N=\{0\}\times(0,1)$,
$\Gamma^D=\partial\Omega\setminus\overline{\Gamma^N}$,
$\varepsilon=10^{-8}$, 
$\bb(x,y)=(-y,x)^\top$, $f=0$, $g=0$, 
\begin{equation*}
  u_b(x,0) = 
  \begin{cases}
  x           & \mathrm{if} \;\; 0       \leq x \leq \frac13 \;, \\
  \frac13 + x & \mathrm{if} \;\; \frac13 < x < \frac23 \;, \\
  1 - x       & \mathrm{if} \;\; \frac23 \leq x \leq 1       \;, \\
  \end{cases}
  \qquad x\in[0,1]\,,
\end{equation*}
and $u_b=0$ elsewhere on $\Gamma^D$. 
\end{example}

The solution contains two interior layers and, due to the small diffusion, the
outflow profile of $u$ at $x=0$ is close to the boundary condition at $y=0$.
Apart from the usual difficulty to approximate the interior layers without
significant smearing and undershoots or overshoots, it is also challenging to
reproduce the piecewise linear parts of the boundary condition correctly.

\begin{figure}[t!]
\begin{center}
\includegraphics[width=0.49\textwidth]{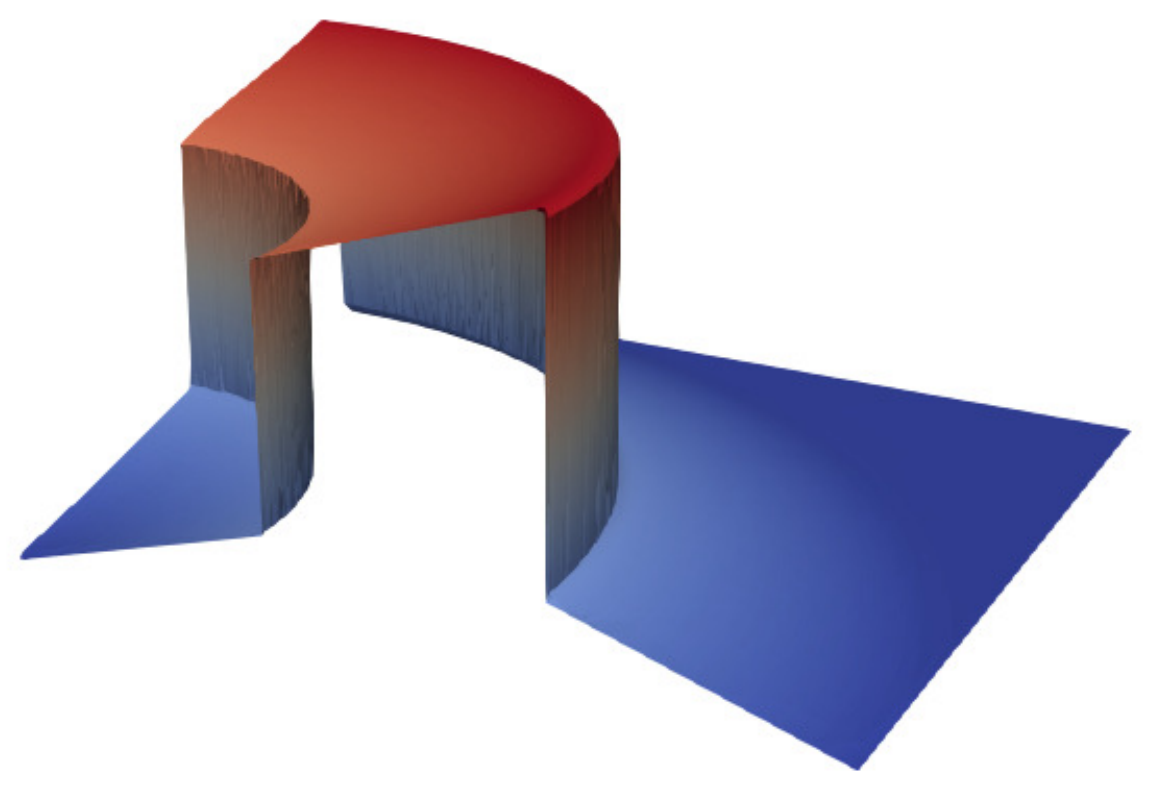}
\includegraphics[width=0.49\textwidth]{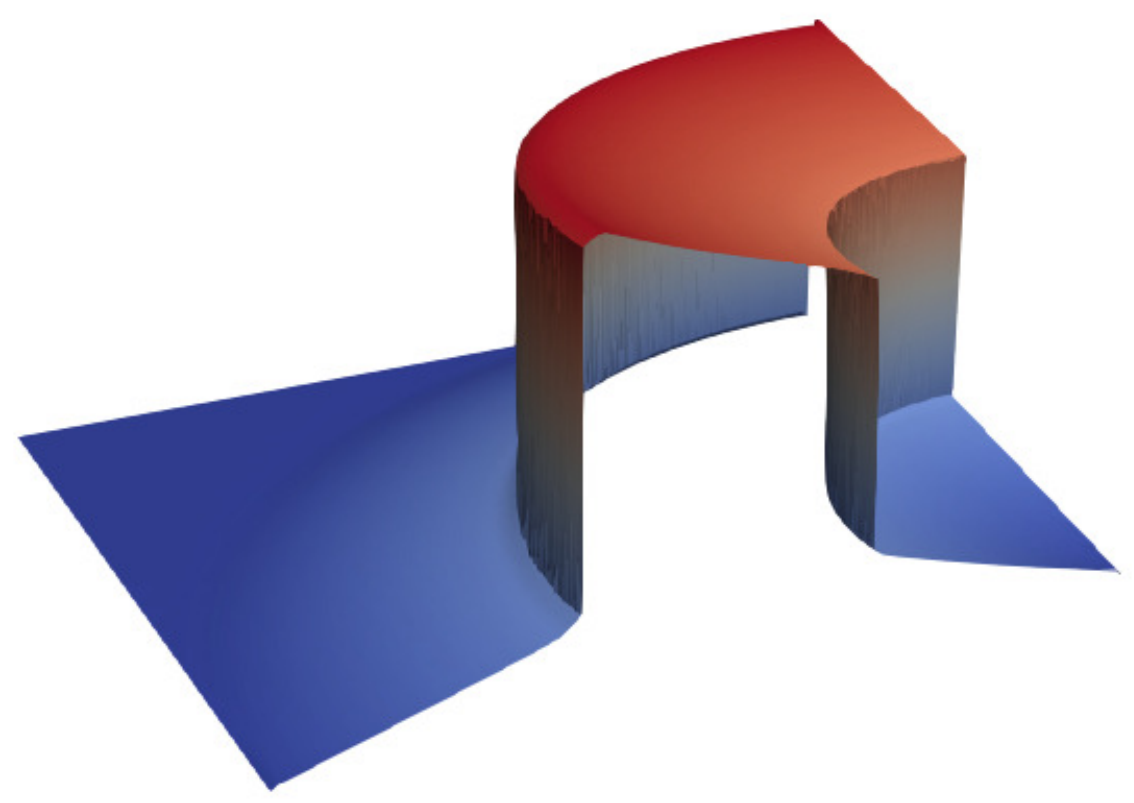}
\end{center}
\caption{Example~\ref{ex2}: two views of the solution of the SUPG method with
$\deltaNN$ obtained on the anisotropic mesh from Fig.~\ref{fig:aniso_meshes}
(right).}
\label{fig:ex2_sol}
\end{figure}

The used triangulation is depicted in Fig.~\ref{fig:aniso_meshes} (right) and
two views of the computed solution are shown in Fig.~\ref{fig:ex2_sol}. The
left-hand graph shows the prescribed inflow profile of the solution, whereas
the right-hand graph shows the computed outflow profile.
\CHANGEF{All features of the solution are captured without spurious
oscillations, and the quantities defined in~(13) listed in
Table~\ref{tab:aniso-ex23} confirm this observation: the overshoot present for
the standard parameters vanishes and the undershoot is reduced to $10^{-4}$.
For this example, the SUPG solution computed with the predicted parameters
exhibited no overshoots relative to the Tabata solution, so that 
no post-processing step was performed.}

Altogether one can state that the approach proposed in
Section~\ref{sec:ml_aniso} represents a promising technique for accurate
solution of convection-dominated problems on anisotropic adapted meshes.

\CHANGEF{

\section{Ablation studies}\label{sec:ablation}

In this section we examine two components of the approach proposed in Section \ref{sec:ml_aniso}:
the choice of the Transformer architecture for the boundary cells and the
post-processing step described at the end of Section \ref{sec:ml_iso}.
Moreover, we compare the ranges of the stabilization parameters given by the
formulas from Section \ref{sec:supg} with those predicted by the neural network
before applying the post-processing step.

\subsection{MLP versus Transformer for the boundary cells}\label{sec:mlp-vs-tf}

As mentioned in Section \ref{sec:ml_aniso}, the MLP architecture used for the interior cells did
not lead to a sufficiently accurate approximation of the target multiplier on the
cells intersecting the outflow Dirichlet boundary, which was the reason for
introducing the tabular Transformer. To quantify this observation, we are going to compare MLP and Transformer on the boundary cells. The interior model and all remaining parts of
the computation were left unchanged.

Table~\ref{tab:mlp-vs-tf} reports the quantities defined in \eqref{eq:min_max_global} for the solutions
of Examples \ref{ex:hemker} and \ref{ex1} obtained on the anisotropic meshes with the two boundary
models, and the corresponding solutions are depicted in Fig \ref{fig:MLPvsTAB}.  With the
MLP, pronounced oscillations appear along the outflow boundary, whereas the
interior layers are approximated equally well by both models, as expected, since
the interior cells are treated by the same network in both cases.

\begin{table}[htbp]
  \centering
  \caption{Examples \ref{ex:hemker} and \ref{ex1}: measures for undershoots and overshoots defined
    by \eqref{eq:min_max_global} for the solutions on the anisotropic meshes obtained with the MLP and
    with the tabular Transformer used on the boundary cells without post processing.}
  \label{tab:mlp-vs-tf}
  \begin{tabular}{|c|c|r|r|}
    \hline
    \textbf{Problem} & \textbf{Quantity} & \textbf{MLP} & \textbf{Transformer}\\ \hline
    Example 1 & $\min_{\mathrm g}$ & 0.0299 &  0.0211\\ \hline
    Example 1 & $\max_{\mathrm g}$ & 0.0358 &  0.0002\\ \hline
    Example 2 & $\min_{\mathrm g}$ & 0.0236 & 0.0051  \\ \hline
    Example 2 & $\max_{\mathrm g}$ & 2.4870 & 0.1542 \\ \hline
  \end{tabular}
\end{table}

\begin{figure}[t!]
\begin{center}
\includegraphics[width=0.45\textwidth]{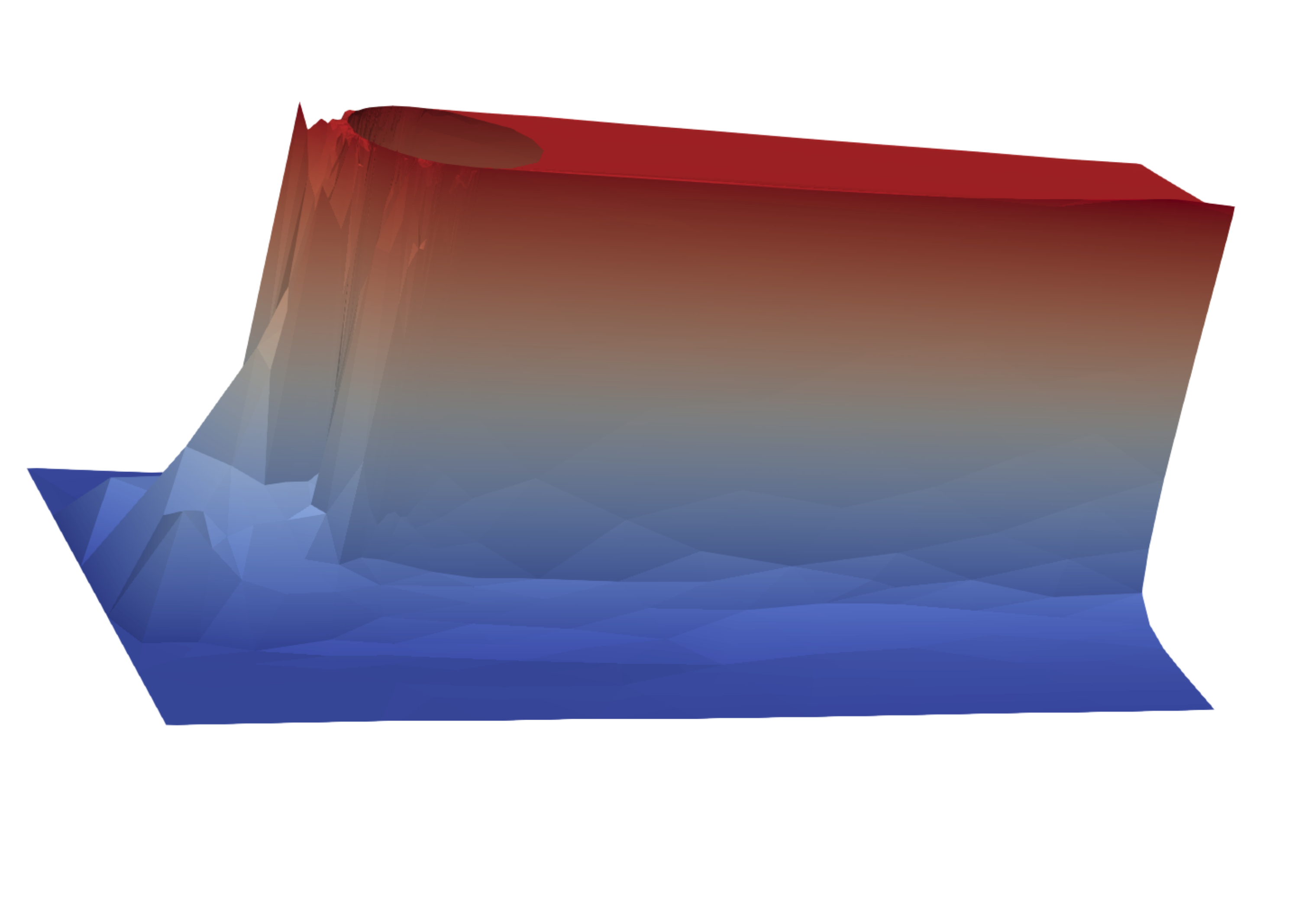}
\hspace*{3mm}
\includegraphics[width=0.50\textwidth]{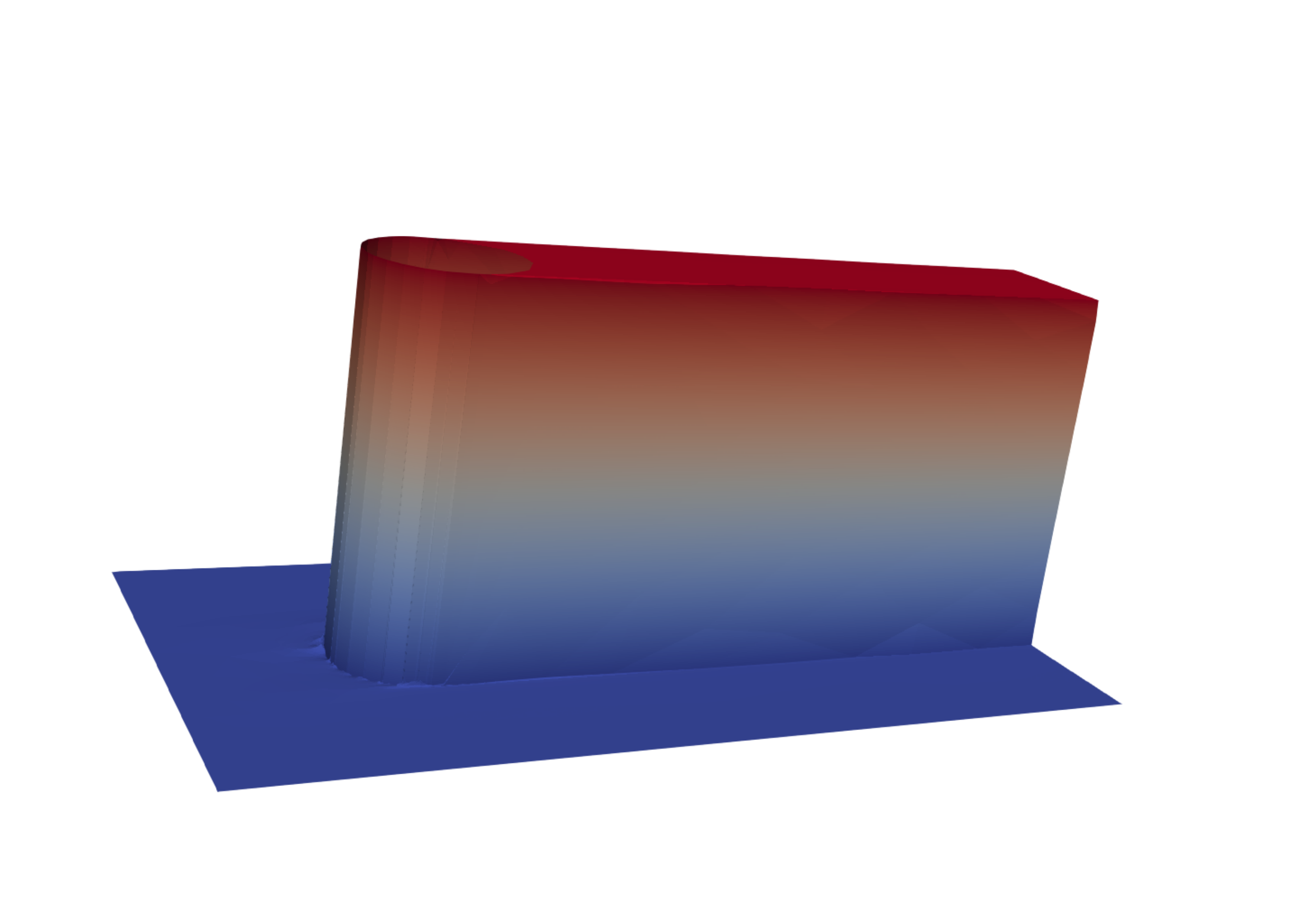}\\[3mm]
\includegraphics[width=0.62\textwidth]{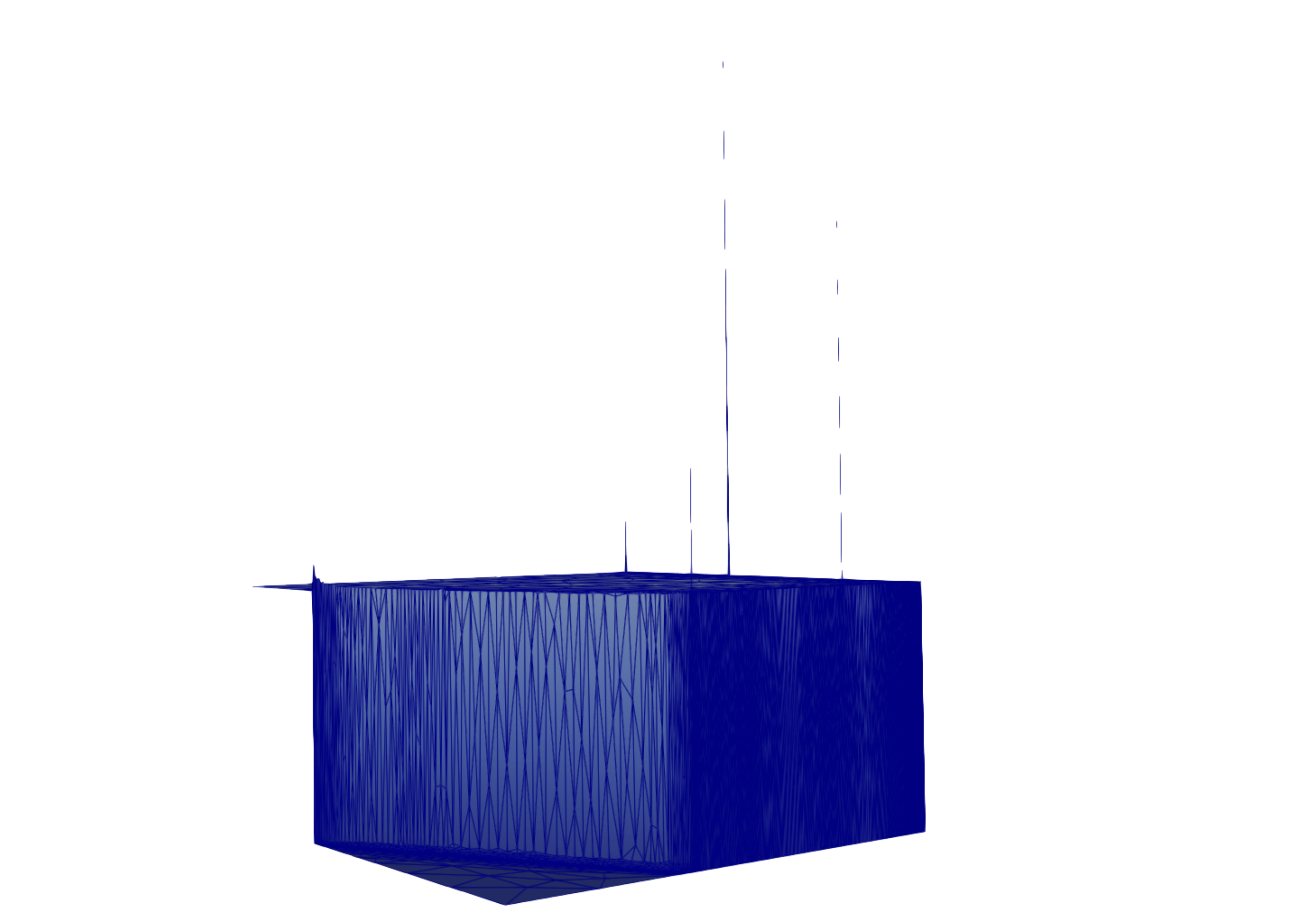}
\hspace*{0mm}
\includegraphics[width=0.33\textwidth]{Figures/sol_supg_opt_8a.pdf}
\end{center}
\caption{Examples~\ref{ex:hemker} and~\ref{ex1}: solutions of the SUPG method with
$\deltaNN$ obtained on the anisotropic meshes, where the boundary cells are
treated by the MLP (left) and by the tabular Transformer (right), for
Example~\ref{ex:hemker} (top) and Example~\ref{ex1} (bottom).}
\label{fig:MLPvsTAB}
\end{figure}

We do not claim a theoretical advantage of self-attention for this problem. The
motivation for the Transformer is the one given in Section \ref{sec:ml_aniso}: the boundary
behaviour of $\delta_K^{\mathrm{opt}}$ depends on more local quantities than in the
interior and the relevant interactions between these quantities are not known a
priori, so that a model learning these interactions during training is a natural
candidate. The comparison above is reported as empirical evidence supporting this
choice.

\subsection{Effect of the post-processing}\label{sec:postproc}

The post-processing step introduced at the end of Section \ref{sec:ml_iso} is applied only on the
cells $K\in\mathscr G_h$ intersecting $\Gamma^{\mathrm{out}}$ and therefore
affects a small number of cells. Its purpose is to remove the overshoots that occasionally
remain along the outflow boundary after the SUPG system has been solved with the
predicted parameters.

Table~\ref{tab:postproc} shows the influence of this step on the solution of
Example \ref{ex1}, and the corresponding solutions are depicted in Fig \ref{fig:NPvsP}.  The limit of five iterations was never reached: two additional solves
were sufficient for Example \ref{ex1}, and for Example \ref{ex2} the solution computed with the
predicted parameters exhibited no overshoots relative to the Tabata solution, so
that the post-processing did not modify the solution at all.

\begin{table}[H]
  \centering
  \caption{Example \ref{ex1}: measures for undershoots and overshoots defined by~\eqref{eq:min_max_global} for
    the solution on the anisotropic mesh computed with the predicted parameters
    without and with the post-processing step.}
  \label{tab:postproc}
  \begin{tabular}{|c|r|r|}
    \hline
    \textbf{Quantity} & \textbf{no post-proc.} & \textbf{post-proc.}\\ \hline
    $\min_{\mathrm g}$ & 0.0051 & 0.0038 \\ \hline
    $\max_{\mathrm g}$ & 0.1542 &  0.0143\\ \hline
  \end{tabular}
\end{table}

\begin{figure}[H]
\begin{center}
\includegraphics[width=0.47\textwidth]{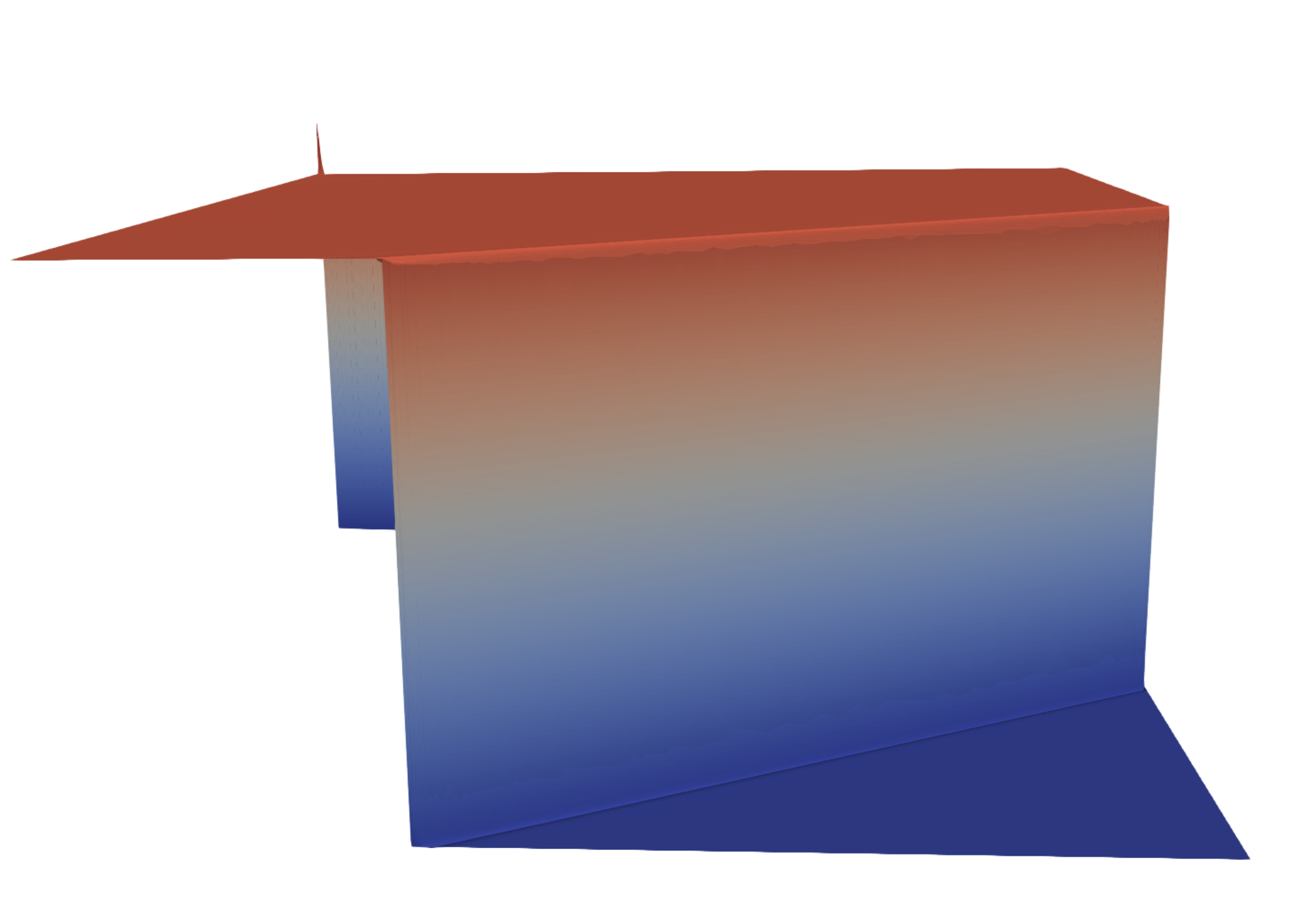}
\hspace*{1mm}
\includegraphics[width=0.5\textwidth]{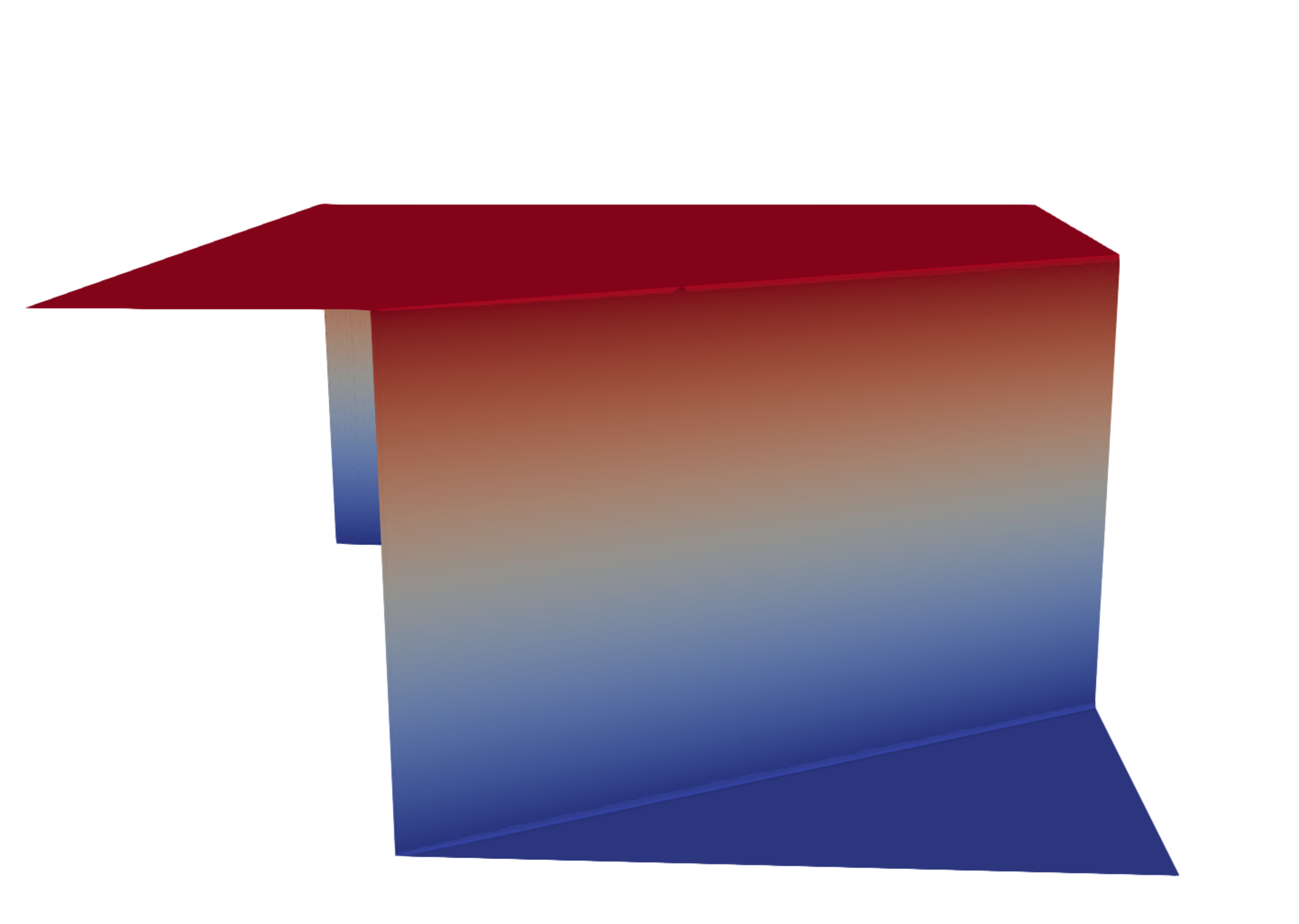}
\end{center}
\caption{Example~\ref{ex1}: solution of the SUPG method with $\deltaNN$ computed
without (left) and with (right) the post-processing step.}
\label{fig:NPvsP}
\end{figure}

\subsection{Ranges of the predicted stabilization parameters}
\label{sec:range}
Table ~\ref{tab:range} compares, for all the computations presented above, the
ranges of the predicted parameters $\delta_K^{\mathrm{NN}}$ with the ranges of
the parameters $\delta_K^{\mathrm{std}}$ and $\delta_K^{\mathrm{bnd}}$, and
presents the ranges of the ratios of $\delta_K^{\mathrm{NN}}$ to the parameters
$\delta_K^{\mathrm{std}}$ and $\delta_K^{\mathrm{bnd}}$. The post-processing
step is not applied. Cells on which $\delta_K^{\mathrm{std}}$ or
$\delta_K^{\mathrm{bnd}}$ vanishes are not included in the ratios. It is
important that the predicted stabilization parameters do not significantly
exceed the parameters $\delta_K^{\mathrm{std}}$ and $\delta_K^{\mathrm{bnd}}$.

\begin{table}[htbp]
  \centering
  \caption{Ranges of $\delta_K^{\mathrm{NN}}$, $\delta_K^{\mathrm{std}}$ and $\delta_K^{\mathrm{bnd}}$,
    and of their ratios, on $K\in\calT_h\setminus\calG_h$ and
    $K\in\calG_h$.}
  \label{tab:range}
  \footnotesize
  \setlength{\tabcolsep}{4pt}
  \begin{tabular}{|l|c|c|c|c|c|c|}
    \hline
    & \multicolumn{3}{c|}{\textbf{$K\in\calT_h\setminus\calG_h$}}
    & \multicolumn{3}{c|}{\textbf{$K\in\calG_h$}}\\
    \textbf{Problem}
      & $\delta_K^{\mathrm{NN}}$ & $\delta_K^{\mathrm{std}}$
      & $\delta_K^{\mathrm{NN}}/\delta_K^{\mathrm{std}}$
      & $\delta_K^{\mathrm{NN}}$ & $\delta_K^{\mathrm{bnd}}$
      & $\delta_K^{\mathrm{NN}}/\delta_K^{\mathrm{bnd}}$\\ \hline
    Ex.\ 1, level 0 & $0.314$--$2.607$ & $0.111$--$0.326$ & $1.61$--$9.67$
                    & $0$--$5.904$ & $0$--$6.092$ & $0.97$--$1.06$\\ \hline
    Ex.\ 1, level 1 & $0.094$--$1.325$ & $0.055$--$0.163$ & $1.50$--$9.69$
                    & $0$--$1.224$ & $0$--$1.259$ & $0.97$--$1.08$\\ \hline
    Ex.\ 1, level 2 & $0.043$--$0.766$ & $0.028$--$0.082$ & $1.50$--$9.70$
                    & $0$--$1.633$ & $0$--$1.685$ & $0.45$--$17.19$\\ \hline
    Ex.\ 1, level 3 & $0.021$--$0.364$ & $0.013$--$0.041$ & $1.50$--$9.69$
                    & $0.011$--$1.248$ & $0.011$--$1.288$ & $0.45$--$18.06$\\ \hline
    Ex.\ 1, level 4 & $0.009$--$0.163$ & $0.006$--$0.020$ & $1.50$--$9.72$
                    & $0.005$--$0.412$ & $0.006$--$0.426$ & $0.45$--$13.35$\\ \hline
    Ex.\ 1, level 5 & $0.004$--$0.080$ & $0.003$--$0.010$ & $1.50$--$9.73$
                    & $0.003$--$0.248$ & $0.003$--$0.256$ & $0.46$--$7.15$\\ \hline
    Ex.\ 1, aniso & $0.0001$--$0.740$ & $0$--$0.094$ & $0.08$--$11.29$
                    & $0.0001$--$0.659$ & $0$--$0.094$ & $0.09$--$9.30$\\ \hline
    Ex.\ 2 & $0.0001$--$0.073$ & $0$--$0.084$ & $0.72$--$32.41$
           & $0$--$0.081$ & $0.0001$--$0.071$ & $0.70$--$3.38$\\ \hline
    Ex.\ 3 & $0.0001$--$0.354$ & $0.0001$--$0.320$ & $0.34$--$3.25$
           & $0.0001$--$0.417$ & $0.0001$--$0.114$ & $0.19$--$3.70$\\ \hline
  \end{tabular}
\end{table}
}

\CHANGEF{

\section{Conclusions and outlook}
\label{sec:conclusions}
We have investigated the machine-learning approach of~\cite{KP26} for computing
SUPG stabilization parameters on unstructured meshes and extended it to meshes
obtained by anisotropic mesh adaptation. On unstructured meshes, the predicted
parameters retain the sharp layer approximation of the standard SUPG method while
reducing the undershoots and overshoots by several orders of magnitude. On
anisotropically adapted meshes, the extended feature sets and the tabular
Transformer model for the boundary cells approximate the reference layer width of
the Hemker problem with an error of $0.0084$ on $5250$ cells, an accuracy
comparable to that obtained on an isotropic triangulation with $760832$ cells.
Both trained models are used for all benchmark problems without retraining, and
none of these problems was used to generate the training data.

The complete pipeline costs about $3$ times a single SUPG solve, essentially
independently of the refinement level, whereas the cost of the a~posteriori
parameter optimization of~\cite{JKS11} grows roughly by an order of magnitude per
refinement level, so that the speedup reaches a factor of about $10^{3}$ on
level~4, see Table~\ref{tab:times-iso}.

Some of the features used in Section~\ref{sec:ml_aniso} are not dimensionless. Replacing these features with dimensionless ones is necessary but is left for future work.

}

\section*{Acknowledgments}
The financial support for the present work was provided by Charles
University project GAUK No. 240325.

\bibliographystyle{alpha}

\newcommand{\etalchar}[1]{$^{#1}$}
\def\cprime{$'$}

\medskip
\medskip

\end{document}